\documentclass[9pt,a4paper,draft]{article}
\usepackage{bbm}
\usepackage{pifont}
\usepackage{bbding}
\usepackage{amsmath}
\usepackage{mathrsfs}
\usepackage{amsfonts}
\usepackage{amssymb}
\usepackage{amsfonts,amssymb,amsmath,indentfirst,amsthm}
\usepackage{epsfig}
\usepackage[numbers,sort&compress]{natbib}
\usepackage{mathtools}
\def\ess~inf{\mathop{\rm ess~inf}}

\numberwithin{equation}{section}

\newenvironment{key words}{\emph{\texttt{Keywords}}\mbox{  }}{ }
\newtheorem{theorem}{Theorem}[section]
\newtheorem{lemma}[theorem]{Lemma}

\theoremstyle{remark}
\newtheorem{remark}[theorem]{\textbf{Remark}}

\theoremstyle{plain}

\makeatletter

\newcommand{\Rmnum}[1]{\expandafter\@slowromancap\romannumeral #1@}
\makeatother

\begin{document}
\setlength{\baselineskip}{1.2\baselineskip}
\title
{\textbf{Pointwise convergence of solutions to the Schr\"{o}dinger equation along a class of curves
} \thanks{This work is supported by Natural Science Foundation of China (No.11601427); Natural Science Basic Research Plan in Shaanxi Province of China (No.2017JQ1009); China Postdoctoral Science Foundation (No.2017M613193); The Fundamental Research Funds for the Central University (No.3102017zy035).}}

\author{Wenjuan Li, \hspace{.2in} Huiju Wang \thanks{Corresponding author's email: huijuwang@mail.nwpu.edu.cn.} \\
\small{Department of Applied Mathematics, Northwestern Polytechnical
University,}\\ \small{ Xi'an, Shaanxi, 710129,  China }}
\date{}
 \maketitle
{\textbf{ Abstract:} In this paper, we obtain pointwise convergence of solutions to the Schrodinger equation along a class of curves in $\mathbb{R}^{2}$ by the polynomial partitioning.}\



\section{\textbf{Introduction}\label{Section 1}}
The solution to the Schr\"{o}dinger equation
\begin{equation}\label{Eq11}
i{u_t} - \Delta u = 0, (x,t) \in {\mathbb{R}^n} \times \mathbb{R},
\end{equation}
with initial datum $u\left( {x,0} \right) = f,$ is formally written as
\[{e^{it\Delta }}f\left( x \right): = \int_{{\mathbb{R}^n}} {{e^{i\left( {x \cdot \xi  + t{{\left| \xi  \right|}^2}} \right)}}\widehat{f}} \left( \xi  \right)d\xi .\]
The problem about finding optimal $s$ for which
\begin{equation}\label{Eq12}
\mathop {\lim }\limits_{t \to 0} {e^{it\Delta }}f\left( x \right) = f(x),a.e.
\end{equation}
whenever $f \in {H^s}\left( {{\mathbb{R}^n}} \right),$ was first considered by Carleson \cite{C}, and extensively studied by Sj\"{o}lin \cite{S} and Vega \cite{V}, who proved independently the convergence for $s > 1/2$ in all dimensions. Dahlberg and Kenig \cite{DK} showed that the convergence does not hold for $s < 1/4$ in any dimension. When $n = 2,$ Du and Li \cite{DL} proved the convergence result for $s > 3/8$ by the polynomial partitioning; Du, Guth and Li \cite{DGL} obtained the sharp result $s > 1/3$ by the polynomial partitioning and ${l^2}$ decoupling method.

By Cho, Lee and Vargas \cite{CLV}, a general generalization of the pointwise convergence problem is to ask a.e. convergence along a wider approach region instead of vertical lines. One of such problems is to consider non-tangential convergence to the initial data, it was shown by Sj\"{o}lin and Sj\"{o}gren \cite{SS} that non-tangential convergence fails for $s \le n/2.$ Another problem is to consider the relation between the degree of the tangency and regularity when $\left( {x,t} \right)$ approaches to $\left( {x,0} \right)$ tangentially.  One of the model problems raised by \cite{CLV} is
 \begin{equation}\label{Eq13}
\mathop {\lim }\limits_{t \to 0} {e^{it\Delta }}f\left( {\gamma \left( {x,t} \right)} \right) = f(x) \:\ a.e.
 \end{equation}
when $n=1$, here the curves $\gamma $ approach $\left( {x,0} \right)$ tangentially to the hyperplane $\left\{ {\left( {x,t} \right):t = 0} \right\}$. Ding and Niu \cite{DN} improved the result of \cite{CLV}, but this problem is still open for $n \ge 2.$

In this paper, we consider this problem when $n = 2$ and
\begin{equation}\label{Eq14}
\gamma \left( {x,t} \right) = x - \sqrt t \mu ,
\end{equation}
where $\mu $ is a unit vector in ${\mathbb{R}^2}.$ The convergence result (\ref{Eq13}) follows from
\begin{theorem}\label{Theorem1.1}
  For $2 \le p \le 3.2,$ if $f \in {H^s}\left({{\mathbb{R}^2}}\right),$ $s > 3/8,$ then there exists a constant $C > 0$ such that
\begin{equation}\label{Eq16}
{\left\| {\mathop {\sup }\limits_{t \in \left( {0,1} \right]} \left|{e^{it\Delta }}f\left( {\gamma \left( {x,t} \right)} \right) \right|} \right\|_{{L^p}\left( {B\left( {0,1} \right)} \right)}} \le C{\left\| f \right\|_{{H^s}}}.
\end{equation}
\end{theorem}

\begin{remark}
When $\mu  = \left( {1,0} \right),$ $s \ge 11/32$ is showed to be necessary if (\ref{Eq16}) holds.
In fact, take
\[\hat{f}\left( \xi  \right): = \psi \left( {\frac{{\xi  - \lambda \mu }}{{{\lambda ^{1/2}}}}} \right),\]
$\psi $ is a non-negative Schwartz function. By rescaling, it follows that
\[\left| {{e^{it\Delta }}f\left( {\gamma \left( {x,t} \right)} \right)} \right| \sim \lambda \]
when $\left| t \right| \le {\lambda ^{ - 1}}$ and $\left| {{\lambda ^{1/2}}\left( {{x_1} + \sqrt t  + 2t\lambda ,{x_2}} \right)} \right| \le C,$ therefore, (\ref{Eq16}) implies that
\[\lambda {\lambda ^{ - 1/2p}} \le {\lambda ^{1/2 + s}}.\]
The desired condition follows from the fact that $\lambda $ can be sufficiently large.
\end{remark}

By Littlewood-Paley Theorem and parabolic rescaling, Theorem \ref{Theorem1.1} can be reduced to
\begin{theorem}\label{Theorem1.2}
 For $2 \le p \le 3.2,$ $\varepsilon  > 0,$ there exists a constant ${C_\varepsilon } > 0$ such that
\begin{equation}\label{Eq17}
{\left\| {\mathop {\sup }\limits_{t \in \left( {0,R} \right]} \left| {e^{it\Delta }}f\left( {\gamma \left( {x,t} \right)} \right) \right|} \right\|_{{L^p}\left( {B\left( {0,R} \right)} \right)}} \le {C_\varepsilon }{R^{\frac{2}{p} - \frac{5}{8} + \varepsilon }}{\left\| f \right\|_{{L^2}}}
\end{equation}
for all $R \ge 1,$ all $f$ with $supp \widehat{f} \subset A\left( 1 \right) = \left\{ {\xi :\left| \xi  \right| = 1} \right\}.$
\end{theorem}

\section{Proof of Theorem \ref{Theorem1.1}}\label{Section 2}
For convenience of the proof, we define a new operator
\begin{equation}\label{Eq15}
{e^{it{\rm H}}}f\left( x \right): = {e^{it\Delta }}f\left( {\gamma \left( {x,t} \right)} \right) = \int_{{\mathbb{R}^2}} {{e^{i\left( {x \cdot \xi  - \sqrt t \mu  \cdot \xi  + t{{\left| \xi  \right|}^2}} \right)}}\widehat{f}} \left( \xi  \right)d\xi.
\end{equation}

\textbf{Proof of Theorem \ref{Theorem1.1}:}  For any $f \in {H^s(\mathbb{R}^{2})},$ we use Littlewood-Paley decomposition,
\begin{align}\label{Eq213}
f = \sum\limits_{k \ge 0} {{f_k},}
  \end{align}
where ${supp \hat{f_{0}}} \subset B\left( {0,1} \right),{supp \hat{f_{k}}} \subset A\left( {{2^k}} \right),k \ge 1.$ If Theorem \ref{Theorem1.2} holds, when $p = 3.2,$ for any $R \ge 1,$ $\hat{g} \in C_c^\infty \left( {{\mathbb{R}^2}} \right)$ with $supp\hat{g} \subset A\left( 1 \right),$ it holds
\[{\left\| {\mathop {\sup }\limits_{t \in \left( {0,R} \right]} \left| {{e^{it{\rm H}}}g} \right|} \right\|_{{L^p}\left( {B\left( {0,R} \right)} \right)}} \le {C_\varepsilon }{R^\varepsilon }{\left\| g \right\|_{{L^2}}},\]
this implies
\begin{align}\label{Eq214}
{\left\| {\mathop {\sup }\limits_{t \in \left( {0,{R^2}} \right]} \left| {{e^{it{\rm H}}}g} \right|} \right\|_{{L^p}\left( {B\left( {0,R} \right)} \right)}} \le {\left\| {\mathop {\sup }\limits_{t \in \left( {0,{R^2}} \right]} \left| {{e^{it{\rm H}}}g} \right|} \right\|_{{L^p}\left( {B\left( {0,{R^2}} \right)} \right)}} \le {C_\varepsilon }{R^{2\varepsilon }}{\left\| g \right\|_{{L^2}}}.
  \end{align}
By parabolic rescaling,
\[\left\{ \begin{array}{l}
 x = Ry, \\
 t = {R^2}s, \\
 \end{array} \right.\]
we have
\begin{align}
 {e^{it{\rm H}}}g\left( x \right) &= \int_{{\mathbb{R}^2}} {{e^{i\left( {x \cdot \xi  - \sqrt t \mu  \cdot \xi  + t{{\left| \xi  \right|}^2}} \right)}}\widehat{g}} \left( \xi  \right)d\xi  = \int_{{\mathbb{R}^2}} {{e^{i\left( {y \cdot R\xi  - \sqrt s \mu  \cdot R\xi  + s{{\left| {R\xi } \right|}^2}} \right)}}\widehat{g}} \left( \xi  \right)d\xi  \nonumber\\
  &= {R^{ - 2}}\int_{{\mathbb{R}^2}} {{e^{i\left( {y \cdot \eta  - \sqrt s \mu  \cdot \eta  + s{{\left| \eta  \right|}^2}} \right)}}\widehat{g}} \left( {\frac{\eta }{R}} \right)d\eta  \nonumber\\
  &= {R^{ - 2}}\int_{{\mathbb{R}^2}} {{e^{i\left( {y \cdot \eta  - \sqrt s \mu  \cdot \eta  + s{{\left| \eta  \right|}^2}} \right)}}\widehat{{g_1}}} \left( \eta  \right)d\eta  \nonumber\\
  &= {R^{ - 2}}{e^{is{\rm H}}}{g_1}\left( y \right), \nonumber
\end{align}
where $\widehat{{g_1}}\left( \eta  \right) := \widehat{g}\left( {\frac{\eta }{R}} \right),$ so that $supp \widehat{{g_1}} \subset A\left( R \right).$ It follows that
\[{\left\| {\mathop {\sup }\limits_{t \in \left( {0,{R^2}} \right]} \left| {{e^{it{\rm H}}}g} \right|} \right\|_{{L^p}\left( {B\left( {0,R} \right)} \right)}} = {R^{2/p - 2}}{\left\| {\mathop {\sup }\limits_{s \in \left( {0,1} \right]} \left| {{e^{is{\rm H}}}{g_1}} \right|} \right\|_{{L^p}\left( {B\left( {0,1} \right)} \right)}}, \hspace{0.2cm} {\left\| g \right\|_{{L^2}}} = {R^{ - 1}}{\left\| {{g_1}} \right\|_{{L^2}}},\]
combining it with (\ref{Eq214}), we have
\begin{align}\label{Eq215}
{\left\| {\mathop {\sup }\limits_{t \in \left( {0,1} \right]} \left| {{e^{it{\rm H}}}g_{1}} \right|} \right\|_{{L^p}\left( {B\left( {0,1} \right)} \right)}} \le {C_\varepsilon }{R^{1 - 2/p + 2\varepsilon }}{\left\| g_{1} \right\|_{{L^2}}}
  \end{align}
for $supp\widehat{g_{1}} \subset A\left( R \right).$ Apply (\ref{Eq215}) to each ${f_k},k \ge 1,$
\begin{align}\label{Eq216}
{\left\| {\mathop {\sup }\limits_{t \in \left( {0,1} \right]} \left| {{e^{it{\rm H}}}{f_k}} \right|} \right\|_{{L^p}\left( {B\left( {0,1} \right)} \right)}} \le {C_\varepsilon }{2^{k{\left( {1 - 2/p + 2\varepsilon } \right)}}}{\left\| {{f_k}} \right\|_{{L^2}}}.
  \end{align}
  And for ${f_0},$ by (\ref{Eq22}) below, it holds
  \begin{align}\label{Eq217}
{\left\| {\mathop {\sup }\limits_{t \in \left( {0,1} \right]} \left| {{e^{it{\rm H}}}{f_0}} \right|} \right\|_{{L^p}\left( {B\left( {0,1} \right)} \right)}} \le {\left\| {{f_0}} \right\|_{{L^2}}}.
  \end{align}
  Combining (\ref{Eq213}), (\ref{Eq216}) and (\ref{Eq217}),
  \begin{align}
 {\left\| {\mathop {\sup }\limits_{t \in \left( {0,1} \right]} \left| {{e^{it{\rm H}}}f} \right|} \right\|_{{L^p}\left( {B\left( {0,1} \right)} \right)}} &\le \sum\limits_{k \ge 0} {{{\left\| {\mathop {\sup }\limits_{t \in \left( {0,1} \right]} \left| {{e^{it{\rm H}}}{f_k}} \right|} \right\|}_{{L^p}\left( {B\left( {0,1} \right)} \right)}}}  \le \sum\limits_{k \ge 0} {{C_\varepsilon }{2^{k{\left( {1 - 2/p + 2\varepsilon } \right)}}}{{\left\| {{f_k}} \right\|}_{{L^2}}}}  \nonumber\\
  &\le \sum\limits_{k \ge 0} {{C_\varepsilon }{2^{k{\left( {1 - 2/p + 2\varepsilon } \right)}}}{2^{ - ks}}{{\left\| {{f}} \right\|}_{{H^s}}}}  \nonumber\\
  &\le C{\left\| f \right\|_{{H^s}}}, \nonumber
  \end{align}
  the last inequality follows from the fact that $s > 3/8$ and $\varepsilon $ can be sufficiently small. Notice that the case $2 \le p < 3.2$ can be easily obtained from the case $p=3.2$ by H\"{o}lder's inequality.

\textbf{Proof of Theorem \ref{Theorem1.2}.} In order to prove (\ref{Eq17}), it suffices to prove that
\begin{equation}\label{Eq21}
{\left\| {\mathop {\sup }\limits_{t \in \left( {0,R} \right]} \left| {{e^{it{\rm H}}}f} \right|} \right\|_{{L^p}\left( {B\left( {0,R} \right)} \right)}} \le {C_\varepsilon }{M^{ - {\varepsilon ^2}}}{R^{\frac{2}{p} - \frac{5}{8} + \varepsilon }}{\left\| f \right\|_{{L^2}}}
\end{equation}
for all $R \ge 1,$ ${\xi _0} \in B\left( {0,1} \right),M \ge 1,$ and any $f$ with $supp\hat{f} \subset B\left( {{\xi _0},{M^{ - 1}}} \right).$

We will prove (\ref{Eq21}) by induction on the physical radius $R$ and frequency radius $1/M$.  So we need to check the base of the induction.

\textbf{Base of the induction.} From
\[\left| {{e^{it{\rm H}}}f} \right| \le {M^{ - 1}}{\left\| f \right\|_{{L^2}}},\]
it is easy to see that
\begin{equation}\label{Eq22}
{\left\| {\mathop {\sup }\limits_{t \in \left( {0,R} \right]} \left| {{e^{it{\rm H}}}f} \right|} \right\|_{{L^p}\left( {B\left( {0,R} \right)} \right)}} \le {M^{ - 1}}{R^{\frac{2}{p}}}{\left\| f \right\|_{{L^2}}}
\end{equation}
for all $R \ge 1,$ $M \ge 1,$ so (\ref{Eq21}) is trivial when $M \ge {R^{10}}.$

When $\sqrt R  \le M \le {R^{10}},$ we adopt wave packets decomposition for $f.$ Let $\varphi $ be a Schwartz function from $\mathbb{R}$ to $\mathbb{R}$, $\hat{\varphi} $ is non-negative and supported in a small neighborhood of the origin, and identically $1$ in another smaller interval. Let $\theta  = \prod\nolimits_{j = 1}^2 {\theta _j} $ denote the rectangle in the frequency space with center $\left( {c\left( {{\theta _1}} \right),c\left( {{\theta _2}} \right)} \right)$ and
\[{\widehat\varphi _\theta }\left( {{\xi _1},{\xi _2}} \right) = \prod\limits_{j = 1}^2 {\frac{1}{{{{\left| {{\theta _j}} \right|}^{1/2}}}}\widehat\varphi } \left( {\frac{{{\xi _j} - c\left( {{\theta _j}} \right)}}{{\left| {{\theta _j}} \right|}}} \right).\]
A rectangle $\nu $ in the physical space is said to be dual to $\theta $ if $\left| {{\theta _j}} \right|\left| {{\nu _j}} \right| = 1,j = 1,2,$ and $\left( {\theta ,\nu } \right)$ is said to be a tile. Let ${\rm T}$ be a collection of all tiles with fixed dimensions and coordinate axes. Define
\[\widehat{{\varphi _{\theta ,\nu }}}\left( \xi  \right) = {e^{ - ic\left( \nu  \right) \cdot \xi }}{\widehat\varphi _\theta }\left( \xi  \right),\]
we have the following representation
\begin{equation}\label{Eq23}
f = \sum\limits_{\left( {\theta ,\nu } \right) \in {\rm T}} {{f_{\theta ,\nu }}}  = \sum\limits_{\left( {\theta ,\nu } \right) \in {\rm T}} {\left\langle {f,{\varphi _{\theta ,\nu }}} \right\rangle } {\varphi _{\theta ,\nu }}.
\end{equation}

We will only use $\left( {\theta ,\nu } \right)$ where $\theta $ is an ${R^{ - 1/2}}$ cube in frequency space and $\nu $ is an ${R^{1/2}}$ cube in physical space. It is clear that
\[\sum\limits_{\left( {\theta ,\nu } \right) \in {\rm T}} {{{\left| {\left\langle {f,{\varphi _{\theta ,\nu }}} \right\rangle } \right|}^2}}  = \left\| f \right\|_{{L^2}}^2.\]
For any Schwartz function $f$ with $supp\hat{f} \subset B\left( {0,1} \right),$ we only need to consider all $\theta's$ that range over $supp\hat{f}$.

Set
\[{\psi _{\theta ,\nu }} = {e^{it{\rm H}}}{\varphi _{\theta ,\nu }},\]
by the representation (\ref{Eq23}), we have
\begin{equation}\label{Eq24}
{e^{it{\rm H}}}f = \sum\limits_{\left( {\theta ,\nu } \right) \in {\rm T}} {{e^{it{\rm H}}}{f_{\theta ,\nu }}}  = \sum\limits_{\left( {\theta ,\nu } \right) \in {\rm T}} {\left\langle {f,{\varphi _{\theta ,\nu }}} \right\rangle } {\psi _{\theta ,\nu }}.
\end{equation}

Next, we consider the localization of ${\psi _{\theta ,\nu }}$ in $B\left( {0,R} \right) \times \left[ {0,R} \right].$ In fact,
\begin{align}
 &{\psi _{\theta ,\nu }}\left( {x,t} \right){\chi _{\left[ {0,R} \right]}}\left( t \right) \nonumber\\
 &= \int_{\mathbb{R}^2} {{e^{i\left( {x \cdot \xi  - \sqrt t \mu  \cdot \xi  + t{{\left| \xi  \right|}^2}} \right)}}{e^{ - ic\left( \nu  \right) \cdot \xi }}{{\widehat\varphi }_\theta }\left( \xi  \right)} d\xi  \times {\chi _{\left[ {0,R} \right]}}\left( t \right) \nonumber\\
  &= \sqrt R \int_{\mathbb{R}^2} {{e^{i\left( {x \cdot \xi  - \sqrt t \mu  \cdot \xi  + t{{\left| \xi  \right|}^2}} \right)}}{e^{ - ic\left( \nu  \right) \cdot \xi }}\prod\limits_{j = 1}^2 {\widehat\varphi } \left( {\frac{{{\xi _j} - c\left( {{\theta _j}} \right)}}{{{R^{ - 1/2}}}}} \right)} d\xi  \times {\chi _{\left[ {0,R} \right]}}\left( t \right) \nonumber\\
  &= \frac{1}{{\sqrt R }}\int_{\mathbb{R}^2} {{e^{i\left( {\left( {x - c\left( \nu  \right)} \right) \cdot \left( {{R^{ - 1/2}}\eta  + c\left( \theta  \right)} \right) - \sqrt t \mu  \cdot \left( {{R^{ - 1/2}}\eta  + c\left( \theta  \right)} \right) + t{{\left| {{R^{ - 1/2}}\eta  + c\left( \theta  \right)} \right|}^2}} \right)}}\prod\limits_{j = 1}^2 {\widehat\varphi } \left( {{\eta _j}} \right)} d\eta  \times {\chi _{\left[ {0,R} \right]}}\left( t \right), \nonumber
\end{align}
the phase function
\[\phi \left( {x,t,\eta } \right) = \left( {x - c\left( \nu  \right)} \right) \cdot \left( {{R^{ - 1/2}}\eta  + c\left( \theta  \right)} \right) - \sqrt t \mu  \cdot \left( {{R^{ - 1/2}}\eta  + c\left( \theta  \right)} \right) + t{\left| {{R^{ - 1/2}}\eta  + c\left( \theta  \right)} \right|^2}.\]
By simple calculation,
\[{\nabla _\eta }\phi \left( {x,t,\eta } \right) = {R^{ - 1/2}}\left( {x - c\left( \nu  \right) + 2tc\left( \theta  \right)} \right) + {R^{ - 1/2}}\sqrt t \mu  + 2{R^{ - 1}}t\eta .\]
It is obvious that in $B\left( {0,R} \right) \times \left[ {0,R} \right],$
\begin{equation}\label{Eq25}
\left| {{\psi _{\theta ,\nu }}\left( {x,t} \right)} \right| \le \frac{1}{{\sqrt R }}{\chi _{{T_{\theta ,\nu }}}}\left( {x,t} \right),
\end{equation}
where ${T_{\theta ,\nu }} := \left\{ {\left( {x,t} \right):0 \le t \le R,\left| {x - c\left( \nu  \right) + 2tc\left( \theta  \right)} \right| \le {R^{\frac{1}{2} + \delta }}} \right\},\delta  \ll \varepsilon ,$ is a tube with direction
\[G\left( \theta  \right) = \left( { - 2c\left( \theta  \right),1} \right).\]

When $M \ge \sqrt R ,$ there is only one possible $\theta $, therefore all tubes are in the same direction. By the definition of $\nu $, $\left| {c\left( {{\nu _1}} \right) - c\left( {{\nu _2}} \right)} \right| \ge {R^{1/2}},{\nu _1} \ne {\nu _2},$ these tubes are also essentially disjoint. What's more, the projection of ${T_{\theta ,\nu }}$ on $x$-plane is contained in an ${R^{1/2}} \times R$ rectangle, denoted by ${S_{\theta ,\nu }},$ by (\ref{Eq25}),
\begin{equation}\label{Eq26}
\left| {{\psi _{\theta ,\nu }}\left( {x,t} \right)} \right| \le \frac{1}{{\sqrt R }}{\chi _{{T_{\theta ,\nu }}}}\left( {x,t} \right) \le \frac{1}{{\sqrt R }}{\chi _{{S_{\theta ,\nu }}}}\left( x \right){\chi _{\left[ {0,R} \right]}}\left( t \right).
\end{equation}
Combining (\ref{Eq24}) and (\ref{Eq26}), we have
\begin{align}
 \left\| {\mathop {\sup }\limits_{t \in \left( {0,R} \right]} \left| {{e^{it{\rm H}}}f} \right|} \right\|_{{L^p}\left( {B\left( {0,R} \right)} \right)}^p &= \int_{B\left( {0,R} \right)}^{} {\mathop {\sup }\limits_{t \in \left( {0,R} \right]} {{\left| {{e^{it{\rm H}}}f} \right|}^p}dx}  \nonumber\\
  &\le \int_{B\left( {0,R} \right)}^{} {\mathop {\sup }\limits_{t \in \left( {0,R} \right]} {{\left| {\sum\limits_{\left( {\theta ,\nu } \right) \in {\rm T}} {\left\langle {f,{\varphi _{\theta ,\nu }}} \right\rangle } {\psi _{\theta ,\nu }}} \right|}^p}dx} \nonumber\\
  &\le \int_{B\left( {0,R} \right)}^{} {\mathop {\sup }\limits_{t \in \left( {0,R} \right]} {{\sum\limits_{\left( {\theta ,\nu } \right) \in {\rm T}} {\left| {\left\langle {f,{\varphi _{\theta ,\nu }}} \right\rangle } \right|} }^p}{{\left| {{\psi _{\theta ,\nu }}} \right|}^p}dx}  \nonumber\\
  &\le {R^{ - p/2}}\int_{B\left( {0,R} \right)}^{} {\mathop {\sup }\limits_{t \in \left( {0,R} \right]} {{\sum\limits_{\left( {\theta ,\nu } \right) \in {\rm T}} {\left| {\left\langle {f,{\varphi _{\theta ,\nu }}} \right\rangle } \right|} }^p}{\chi _{{S_{\theta ,\nu }}}}\left( x \right){\chi _{\left[ {0,R} \right]}}\left( t \right)dx}  \nonumber\\
  &\le {R^{\frac{{3 - p}}{2} + O\left( \delta  \right)}}\left\| f \right\|_{{L^2}}^p, \nonumber
\end{align}
from which (\ref{Eq21}) follows.

Therefore, we only need to consider the case $M \le \sqrt R .$  On the other hand, when $R \le C$
 for some constant $C > 0,$ the result is true by (\ref{Eq22}). So we can assume that $R$ is sufficiently large. This completes the base of our induction. Now we are ready to prove Theorem \ref{Theorem1.2}.

 Choose non-negative Schwartz functions ${\psi _1}\left( t \right)$
 and ${\psi _2}\left( t \right)$, such that ${\psi _1}\left( t \right)$ is supported in a sufficiently small neighborhood of $\left[ {0,{R^{\varepsilon  - 1}}} \right],$ and identically $1$ on $\left[ {0,{R^{\varepsilon  - 1}}} \right],$ ${\psi _2}\left( t \right)$ is supported in a sufficiently small neighborhood of $\left[ {{R^{\varepsilon  - 1}},1} \right],$ and identically $1$ on $\left[ {{R^{\varepsilon  - 1}},1} \right].$ We have
\begin{align}\label{Eq27}
 {\left\| {\mathop {\sup }\limits_{t \in \left( {0,R} \right]} \left| {{e^{it{\rm H}}}f} \right|} \right\|_{{L^p}\left( {B\left( {0,R} \right)} \right)}} &\le {\left\| {\mathop {\sup }\limits_{t \in \left( {0,R} \right]} \left| {{e^{it{\rm H}}}f} \right|{\psi _1}\left( {\frac{t}{R}} \right) + \mathop {\sup }\limits_{t \in \left( {0,R} \right]} \left| {{e^{it{\rm H}}}f} \right|{\psi _2}\left( {\frac{t}{R}} \right)} \right\|_{{L^p}\left( {B\left( {0,R} \right)} \right)}} \nonumber\\
  &\le {\left\| {\mathop {\sup }\limits_{t \in \left( {0,R} \right]} \left| {{e^{it{\rm H}}}f} \right|{\psi _1}\left( {\frac{t}{R}} \right)} \right\|_{{L^p}\left( {B\left( {0,R} \right)} \right)}} + {\left\| {\mathop {\sup }\limits_{t \in \left( {0,R} \right]} \left| {{e^{it{\rm H}}}f} \right|{\psi _2}\left( {\frac{t}{R}} \right)} \right\|_{{L^p}\left( {B\left( {0,R} \right)} \right)}} \nonumber\\
 & := {I_1} + {I_2}.
\end{align}

If ${I_1}$ dominates, $t$ is localized in a sufficiently small neighborhood of $\left[ {0,{R^{\varepsilon}}} \right],$ the oscillatory integral
\[{e^{it{\rm H}}}f = \int_{\mathbb{R}^2} {{e^{i\left( {x \cdot \xi  - \sqrt t \mu  \cdot \xi  + t{{\left| \xi  \right|}^2}} \right)}}\hat{f}} \left( \xi  \right)d\xi \]
shows that $\left| {{e^{it{\rm H}}}f} \right|{\psi _1}\left( {\frac{t}{R}} \right)$ is essentially supported in $B\left( {0,{R^{1 - \varepsilon }}} \right) \times \left[ {0,{R^{\varepsilon }}} \right].$ Therefore,
\begin{align}\label{Eq28}
 {I_1} &\le {\left\| {\mathop {\sup }\limits_{t \in \left( {0,{R^{1 - \varepsilon }}} \right]} \left| {{e^{it{\rm H}}}f} \right|{\psi _1}\left( {\frac{t}{R}} \right)} \right\|_{{L^p}\left( {B\left( {0,{R^{1 - \varepsilon }}} \right)} \right)}} \le 2{C_\varepsilon }{M^{ - {\varepsilon ^2}}}{R^{{{\left( {1 - \varepsilon } \right)}{\left( {\frac{2}{p} - \frac{5}{8} + \varepsilon } \right)}}}}{\left\| f \right\|_{{L^2}}} \nonumber\\
  &\le {R^{ - {\varepsilon ^2}}}{C_\varepsilon }{M^{ - {\varepsilon ^2}}}{R^{\frac{2}{p} - \frac{5}{8} + \varepsilon }}{\left\| f \right\|_{{L^2}}},
 \end{align}
since $R$ is sufficiently large, then  (\ref{Eq27}) and (\ref{Eq28}) finished the induction.

We consider the case when ${I_2}$ dominates.  Let $K$ be a large parameter such that $K \ll {R^\delta },$ we decompose $B\left( {0,R} \right)$ into balls ${B_K}$ of  radius $K$, and interval $\left[ {0,R} \right]$ into intervals $I_K^j$ of length $K$. We write
\begin{equation}\label{Eq29}
\left\| {\mathop {\sup }\limits_{t \in \left( {0,R} \right]} \left| {{e^{it{\rm H}}}f} \right|{\psi _2}\left( {\frac{t}{R}} \right)} \right\|_{{L^p}\left( {B\left( {0,R} \right)} \right)}^p = \sum\limits_{{B_K} \subset B\left( {0,R} \right)} {\int_{{B_K}} {\mathop {\sup }\limits_{I_K^j \subset \left[ {0,R} \right]} } } \mathop {\sup }\limits_{t \in I_K^j} {\left| {{e^{it{\rm H}}}f} \right|^p}{\psi _2}{\left( {\frac{t}{R}} \right)^p}dx.
\end{equation}
We divide $B\left( {{\xi _0},{M^{ - 1}}} \right)$ into balls $\tau $ of radius ${\left( {KM} \right)^{ - 1}},$ $f = \sum\nolimits_\tau  {{f_\tau }} ,\widehat{{f_\tau }} = {\left. {\hat{f}} \right|_\tau }.$ For each ${B_K} \times I_K^j$ and a parameter $A \in \mathbb{Z}^{+}$, we choose 1-dimensional sub-spaces $V_1^{0},V_2^{0},...,V_A^{0}$ such that
\begin{equation}\label{Eq210}
{\mu _{{e^{it{\rm H}}}f{\psi _2}\left( {\frac{t}{R}} \right)}}\left( {{B_K} \times I_K^j} \right) := \mathop {\min }\limits_{{V_1},{V_2},...,{V_A}} \left( {\mathop {\max }\limits_{\tau  \notin {V_\alpha },\alpha  = 1,2,...,A} \int_{{B_K} \times I_K^j} {{{\left| {{e^{it{\rm H}}}{f_\tau }} \right|}^p}{\psi _2}{{\left( {\frac{t}{R}} \right)}^p}dxdt} } \right)
\end{equation}
achieves the minimum. We say that $\tau  \in {V_\alpha }$ if
\[\mathop {\inf }\limits_{\xi  \in \tau } Angle\left( {\frac{{\left( { - 2\xi ,1} \right)}}{{\left| {\left( { - 2\xi ,1} \right)} \right|}},{V_\alpha }} \right) \le {\left( {KM} \right)^{ - 1}}.\]
Then from (\ref{Eq29}),
\begin{align}
 &\left\| {\mathop {\sup }\limits_{t \in \left( {0,R} \right]} \left| {{e^{it{\rm H}}}f} \right|{\psi _2}\left( {\frac{t}{R}} \right)} \right\|_{{L^p}\left( {B\left( {0,R} \right)} \right)}^p \nonumber\\
 &= \sum\limits_{{B_K} \subset B\left( {0,R} \right)} {\int_{{B_K}} {\mathop {\sup }\limits_{I_K^j \subset \left[ {0,R} \right]} } } \mathop {\sup }\limits_{t \in I_K^j} {\left| {\sum\nolimits_\tau  {{e^{it{\rm H}}}{f_\tau }} } \right|^p}{\psi _2}{\left( {\frac{t}{R}} \right)^p}dx \nonumber\\
  &\le \sum\limits_{{B_K} \subset B\left( {0,R} \right)} {\int_{{B_K}} {\mathop {\sup }\limits_{I_K^j \subset \left[ {0,R} \right]} } } \mathop {\sup }\limits_{t \in I_K^j} {\left| {\sum\nolimits_{\tau  \notin V_\alpha^{0},\alpha  = 1,2,...,A} {{e^{it{\rm H}}}{f_\tau }} } \right|^p}{\psi _2}{\left( {\frac{t}{R}} \right)^p}dx \nonumber\\
  &+ \sum\limits_{{B_K} \subset B\left( {0,R} \right)} {\int_{{B_K}} {\mathop {\sup }\limits_{I_K^j \subset \left[ {0,R} \right]} } } \mathop {\sup }\limits_{t \in I_K^j} {\left| {\sum\nolimits_{\tau  \in \;some\;V_\alpha^{0},\alpha  = 1,2,...,A} {{e^{it{\rm H}}}{f_\tau }} } \right|^p}{\psi _2}{\left( {\frac{t}{R}} \right)^p}dx \nonumber\\
 & := {I_3} + {I_4}. \nonumber
\end{align}

 If ${I_3}$ dominates, we have
\begin{align}\label{Eq211}
 {I_3} &\le \sum\limits_{{B_K} \subset B\left( {0,R} \right)} {\int_{{B_K}} {\mathop {\sup }\limits_{I_K^j \subset \left[ {0,R} \right]} } } \mathop {\sup }\limits_{t \in I_K^j} {\left| {\sum\nolimits_{\tau  \notin V_\alpha^{0},\alpha  = 1,2,...,A} {{e^{it{\rm H}}}{f_\tau }} } \right|^p}{\psi _2}{\left( {\frac{t}{R}} \right)^p}dx \nonumber\\
  &\le {K^{{\rm O}\left( 1 \right)}}\sum\limits_{{B_K} \subset B\left( {0,R} \right)} {\int_{{B_K}} {\mathop {\sup }\limits_{I_K^j \subset \left[ {0,R} \right]} } } \mathop {\sup }\limits_{t \in I_K^j} \mathop {\max }\limits_{\tau  \notin V_\alpha^{0},\alpha  = 1,2,...,A} {\left| {{e^{it{\rm H}}}{f_\tau }} \right|^p}{\psi _2}{\left( {\frac{t}{R}} \right)^p}dx \nonumber\\
  &= {K^{{\rm O}\left( 1 \right)}}\sum\limits_{{B_K} \subset B\left( {0,R} \right)} {\int_{{B_K}} {\mathop {\sup }\limits_{I_K^j \subset \left[ {0,R} \right]} } } \mathop {\sup }\limits_{t \in I_K^j} \mathop {\max }\limits_{\tau  \notin V_\alpha^{0},\alpha  = 1,2,...,A} c_K^jdx \nonumber\\
  &= {K^{{\rm O}\left( 1 \right)}}\sum\limits_{{B_K} \subset B\left( {0,R} \right)} {\mathop {\sup }\limits_{I_K^j \subset \left[ {0,R} \right]} } \mathop {\min }\limits_{{V_1},{V_2},...,{V_A}} \left( {\mathop {\max }\limits_{\tau  \notin {V_\alpha },\alpha  = 1,2,...,A} \int_{{B_K} \times I_K^j} {{{\left| {{e^{it{\rm H}}}{f_\tau }} \right|}^p}{\psi _2}{{\left( {\frac{t}{R}} \right)}^p}dxdt} } \right),
  \end{align}
where we used the fact that ${\left| {{e^{it{\rm H}}}{f_\tau }} \right|^p}{\psi _2}\left( {\frac{t}{R}} \right)$  is essentially constant $c_{K}^{j}$ on ${B_K} \times I_K^j.$ Denote
\begin{align}\label{Eq211}
&\left\| {{e^{it{\rm H}}}f{\psi _2}\left( {\frac{t}{R}} \right)} \right\|_{BL_{k,A}^p{L^\infty }\left( {B\left( {0,R} \right) \times \left[ {0,R} \right]} \right)}^p \nonumber\\
&: = \sum\limits_{{B_K} \subset B\left( {0,R} \right)} {\mathop {\sup }\limits_{I_K^j \subset \left[ {0,R} \right]} } \mathop {\min }\limits_{{V_1},{V_2},...,{V_A}} \left( {\mathop {\max }\limits_{\tau  \notin {V_\alpha },\alpha  = 1,2,...,A} \int_{{B_K} \times I_K^j} {{{\left| {{e^{it{\rm H}}}{f_\tau }} \right|}^p}{\psi _2}{{\left( {\frac{t}{R}} \right)}^p}dxdt} } \right),\nonumber
  \end{align}
by Theorem \ref{Theorem 2.1} below, we have
\[{I_3} \le {K^{{\rm O}\left( 1 \right)}}{R^{\frac{p}{2}{\varepsilon ^2}}}{\left[ {C\left( {K,\frac{\varepsilon }{2}} \right){M^{ - {\varepsilon ^2}}}{R^{\frac{2}{p} - \frac{5}{8} + \frac{\varepsilon }{2}}}{{\left\| f \right\|}_{{L^2}}}} \right]^p},\]
(\ref{Eq21}) follows from the fact that $R$ is sufficiently large. If $I_{4}$ dominates, we have
\begin{align}
 {I_4} &\le \sum\limits_{{B_K} \subset B\left( {0,R} \right)} {\int_{{B_K}} {\mathop {\sup }\limits_{I_K^j \subset \left[ {0,R} \right]} } } \mathop {\sup }\limits_{t \in I_K^j} {\left( {\sum\nolimits_{\tau  \in \;some\;V_\alpha^{0},\alpha  = 1,2,...,A} {\left| {{e^{it{\rm H}}}{f_\tau }} \right|} } \right)^p}dx \nonumber\\
  &\le \sum\limits_{{B_K} \subset B\left( {0,R} \right)} {\int_{{B_K}} {\mathop {\sup }\limits_{I_K^j \subset \left[ {0,R} \right]} } } \mathop {\sup }\limits_{t \in I_K^j} {\left( {\sum\nolimits_{\alpha  = 1}^A {\sum\nolimits_{\tau  \in V_\alpha^{0}}^{} {\left| {{e^{it{\rm H}}}{f_\tau }} \right|} } } \right)^p}dx \nonumber\\
  &\le {A^{p - 1}}\sum\limits_{{B_K} \subset B\left( {0,R} \right)} {\int_{{B_K}} {\mathop {\sup }\limits_{I_K^j \subset \left[ {0,R} \right]} } } \mathop {\sup }\limits_{t \in I_K^j} \sum\nolimits_{\alpha  = 1}^A {\sum\nolimits_{\tau  \in V_\alpha^{0}}^{} {{{\left| {{e^{it{\rm H}}}{f_\tau }} \right|}^p}} } dx \nonumber\\
  &\le O\left( 1 \right){A^{p - 1}}A\sum\nolimits_\tau ^{} {\sum\limits_{{B_K} \subset B\left( {0,R} \right)} {\int_{{B_K}} {\mathop {\sup }\limits_{I_K^j \subset \left[ {0,R} \right]} } } \mathop {\sup }\limits_{t \in I_K^j} {{\left| {{e^{it{\rm H}}}{f_\tau }} \right|}^p}dx}  \nonumber\\
  &\le O\left( 1 \right){A^p}{K^{ - {\varepsilon ^2}p}}\sum\nolimits_\tau ^{} {{{\left( {{C_\varepsilon }{M^{ - {\varepsilon ^2}}}{R^{\frac{2}{p} - \frac{5}{8} + \varepsilon }}{{\left\| {{f_\tau }} \right\|}_{{L^2}}}} \right)}^p}}  \nonumber\\
  &\le O\left( 1 \right){A^p}{K^{ - {\varepsilon ^2}p}}\left( {{C_\varepsilon }{M^{ - {\varepsilon ^2}}}{R^{\frac{2}{p} - \frac{5}{8} + \varepsilon }}{{\left\| {{f }} \right\|}_{{L^2}}}} \right)^{p},\nonumber
  \end{align}
choose $K$ sufficiently large such that $A{K^{ - {\varepsilon ^2}}} \ll 1.$

In the proof of Theorem \ref{Theorem1.2}, we used
\begin{theorem}\label{Theorem 2.1}
 For $2 \le p \le 3.2$ and $k=2,$ for any $\varepsilon  > 0,$ there exist positive constants $A = A\left( \varepsilon  \right)$ and $C\left( {K,\varepsilon } \right)$ such that
\begin{equation}\label{Eq212}
\left\| {{e^{it{\rm H}}}f{\psi _2}\left( {\frac{t}{R}} \right)} \right\|_{BL_{k,A}^p{L^\infty }\left( {B\left( {0,R} \right) \times \left[ {0,R} \right]} \right)}^{} \le C\left( {K,\varepsilon } \right){R^{\frac{2}{p} - \frac{5}{8} + \varepsilon }}{\left\| f \right\|_{{L^2}}},
\end{equation}
for all $R \ge 1,$ ${\xi _0} \in B\left( {0,1} \right),M \ge 1,$ all $f$ with $supp\hat{f} \subset B\left( {{\xi _0},{M^{ - 1}}} \right).$
\end{theorem}

We will prove Theorem \ref{Theorem 2.1} from Section 3 to Section 7.

\section{Preliminaries for the proof of Theorem \ref{Theorem 2.1}}\label{Section 3}
For any subset $U \subset B\left( {0,R} \right) \times \left[ {0,R} \right],$ we define
\[\left\| {{e^{it{\rm H}}}f{\psi _2}\left( {\frac{t}{R}} \right)} \right\|_{BL_{k,A}^p{L^\infty }\left( U \right)}^{} := {\left( {\sum\limits_{{B_K} \subset B\left( {0,R} \right)} {\mathop {\sup }\limits_{I_K^j \subset \left[ {0,R} \right]} \frac{{\left| {U \cap \left( {{B_K} \times I_K^j} \right)} \right|}}{{\left| {{B_K} \times I_K^j} \right|}}{\mu _{{e^{it{\rm H}}}{f}{\psi _2}\left( {\frac{t}{R}} \right)}}\left( {{B_K} \times I_K^j} \right)} } \right)^{1/p}},\]
which can be approximated by
\begin{align}
&\left\| {{e^{it{\rm H}}}f{\psi _2}\left( {\frac{t}{R}} \right)} \right\|_{BL_{k,A}^p{L^q}\left( U \right)} \nonumber\\
:&= {\left( {{{\sum\limits_{{B_K} \subset B\left( {0,R} \right)} {\left[ {\sum\limits_{I_K^j \subset \left[ {0,R} \right]}^{} {{{\left( {\frac{{\left| {U \cap \left( {{B_K} \times I_K^j} \right)} \right|}}{{\left| {{B_K} \times I_K^j} \right|}}{\mu _{{e^{it{\rm H}}}{f}{\psi _2}\left( {\frac{t}{R}} \right)}}\left( {{B_K} \times I_K^j} \right)} \right)}^q}} } \right]} }^{1/q}}} \right)^{1/p}},\nonumber
\end{align}
i.e.,
\[\left\| {{e^{it{\rm H}}}f{\psi _2}\left( {\frac{t}{R}} \right)} \right\|_{BL_{k,A}^p{L^\infty }\left( U \right)}^{} = \mathop {\lim }\limits_{q \to  + \infty } \left\| {{e^{it{\rm H}}}f{\psi _2}\left( {\frac{t}{R}} \right)} \right\|_{BL_{k,A}^p{L^q }\left( U \right)}^{},\]
which implies that Theorem \ref{Theorem 2.1} can be turned to prove Theorem \ref{Theorem3.1}.
\begin{theorem}\label{Theorem3.1}
 For $2 \le p \le 3.2$ and $k=2$, for any $\varepsilon  > 0,$ $1 \le q <  + \infty,$ there exist positive constants $A = A\left( \varepsilon  \right)$ and $C\left( {K,\varepsilon } \right)$ such that
\[\left\| {{e^{it{\rm H}}}f{\psi _2}\left( {\frac{t}{R}} \right)} \right\|_{BL_{k,A}^p{L^q}\left( {B\left( {0,R} \right) \times \left[ {0,R} \right]} \right)}^{} \le C\left( {K,\varepsilon } \right){R^{\frac{1}{{qp}}}}{R^{\frac{2}{p} - \frac{5}{8} + \varepsilon }}{\left\| f \right\|_{{L^2}}},\]
for all $R \ge 1,$ ${\xi _0} \in B\left( {0,1} \right),M \ge 1,$ all $f$ with $supp \widehat{f} \subset B\left( {{\xi _0},{M^{ - 1}}} \right).$
\end{theorem}
Instead of Theorem \ref{Theorem3.1}, we will prove Theorem \ref{Theorem3.2} below.
 \begin{theorem}\label{Theorem3.2}
For $2 \le p \le 3.2$ and $k=2$, for any $\varepsilon  > 0,$ $1 \le q <  + \infty,$ there exist positive constants $\overline {A}  = \overline {A} \left( \varepsilon  \right)$ and $C\left( {K,\varepsilon } \right)$ such that
\begin{equation}\label{Eq31}
\left\| {{e^{it{\rm H}}}f{\psi _2}\left( {\frac{t}{R}} \right)} \right\|_{BL_{k,A}^p{L^q}\left( {B\left( {0,{R^{'}}} \right) \times [0,R] } \right)}^{} \le C\left( {K,\varepsilon } \right){R^{\delta \left( {\log \overline A  - \log A} \right)}}{R^{\frac{1}{{qp}}}}{\left( {{R^{'}}} \right)^{\frac{2}{p} - \frac{5}{8} + \varepsilon }}{\left\| f \right\|_{{L^2}}},
\end{equation}
for any fixed $R \ge 1,$ all $1 \le {R^{'}} \le R,$ $1 \le A \le \overline A ,$ ${\xi _0} \in B\left( {0,1} \right),M \ge 1,$ all $f$ with $supp \widehat{f} \subset B\left( {{\xi _0},{M^{ - 1}}} \right).$
\end{theorem}
We will prove Theorem \ref{Theorem3.2} by induction on $R^{'}$ and $A$, we will check the base of the induction.

\textbf{   The base of the induction.}  Given $R > 1,$ for any $1 \le {R^{'}} \le R,$ it is easy to see
\begin{equation}\label{Eq32}
\left\| {{e^{it{\rm H}}}f{\psi _2}\left( {\frac{t}{R}} \right)} \right\|_{BL_{k,A}^p{L^q}\left( {B\left( {0,{R^{'}}} \right) \times [0,R]} \right)}^{} \le C\left( K \right){R^{1/pq}}{\left( {{R^{'}}} \right)^{\frac{2}{p}}}\left\| f \right\|_{{L^2}}^{}.
\end{equation}

(1) When ${R^{'}}$ is controlled by some constant $C$, then (\ref{Eq31}) holds.

(2) When $A = 1,$ then (\ref{Eq31}) holds even though $A$ does not appear in the right side of (\ref{Eq32}). In fact, we choose $\overline {A} $ such that $\delta \log \overline {A}  = 100$, therefore
 \begin{align}
 \left\| {{e^{it{\rm H}}}f{\psi _2}\left( {\frac{t}{R}} \right)} \right\|_{BL_{k,A}^p{L^q}\left( {B\left( {0,{R^{'}}} \right) \times [0,R]} \right)}^{} &\le C\left( K \right){R^{1/pq}}{\left( {{R^{'}}} \right)^{\frac{5}{8} - \varepsilon }}{\left( {{R^{'}}} \right)^{\frac{2}{p} - \frac{5}{8} + \varepsilon }}\left\| f \right\|_{{L^2}}^{} \nonumber\\
   &\le C\left( K \right){R^{1/pq}}{R^{100}}{\left( {{R^{'}}} \right)^{\frac{2}{p} - \frac{5}{8} + \varepsilon }}\left\| f \right\|_{{L^2}}^{} \nonumber\\
  &= C\left( K \right){R^{1/pq}}{R^{\delta \left( {\log \overline {A}  - \log A} \right)}}{\left( {{R^{'}}} \right)^{\frac{2}{p} - \frac{5}{8} + \varepsilon }}\left\| f \right\|_{{L^2}}^{}, \nonumber
  \end{align}
  this completes the base of the induction. What's more, by the analysis in Section 2, we only need to consider the case $KM \le {R^{1/2}}.$

  In order to prove Theorem \ref{Theorem3.2}, we need some basic inequalities:
  \begin{lemma}\label{Lemma3.3}
(1) If ${U_1}$ and ${U_2}$ are two subsets of $B\left( {0,R} \right) \times \left[ {0,R} \right],$
 then for $1 \le q <  + \infty ,$
\[\left\| {{e^{it{\rm H}}}f{\psi _2}\left( {\frac{t}{R}} \right)} \right\|_{BL_{k,A}^p{L^q}\left( {{U_1} \cup {U_2}} \right)}^p \le \left\| {{e^{it{\rm H}}}f{\psi _2}\left( {\frac{t}{R}} \right)} \right\|_{BL_{k,A}^p{L^q}\left( {{U_1}} \right)}^p + \left\| {{e^{it{\rm H}}}f{\psi _2}\left( {\frac{t}{R}} \right)} \right\|_{BL_{k,A}^p{L^q}\left( {{U_2}} \right)}^p.\]
(2) Given non-negative integers $A,{A_1},{A_2},$ $A = {A_1} + {A_2},$ then for $1 \le q <  + \infty ,$
\[\left\| {{e^{it{\rm H}}}\left( {f + g} \right){\psi _2}\left( {\frac{t}{R}} \right)} \right\|_{BL_{k,A}^p{L^q}\left( U \right)}^p \le {C_p}\left( {\left\| {{e^{it{\rm H}}}f{\psi _2}\left( {\frac{t}{R}} \right)} \right\|_{BL_{k,{A_1}}^p{L^q}\left( U \right)}^p + \left\| {{e^{it{\rm H}}}g{\psi _2}\left( {\frac{t}{R}} \right)} \right\|_{BL_{k,{A_2}}^p{L^q}\left( U \right)}^p} \right).\]
(3) If $1 \le p \le r,$ $U \subset {S_U} \times {I_U} \subset B(0,R) \times [0,R],$ where $S_{U}$ and $I_{U}$ are subsets paralleled to the $x$-plane and $t$-axe respectively, then for $1 \le q <  + \infty ,$
\[\left\| {{e^{it{\rm H}}}f{\psi _2}\left( {\frac{t}{R}} \right)} \right\|_{BL_{k,A}^p{L^q}\left( U \right)}^{} \le {C_K}{\left( {\left| {{S_U}} \right|{{\left| {{I_U}} \right|}^{1/q}}} \right)^{\left( {\frac{1}{p} - \frac{1}{r}} \right)}}\left\| {{e^{it{\rm H}}}f{\psi _2}\left( {\frac{t}{R}} \right)} \right\|_{BL_{k,A}^r{L^q}\left( U \right)}^{}.\]
\end{lemma}

The proof of Lemma \ref{Lemma3.3} is very similar to Lemma 3.1, \cite{DL}, So we omit the proof here.

  \section{Polynomial partitioning}\label{Section 4}
The main tool we will use is polynomial partitioning.
\begin{lemma}\label{Lemma4.1}
Suppose ${f_1},{f_2},...,{f_N}$ are functions defined on $\mathbb{R}^{n}$ with $supp \widehat{{f_j}} \subset {B^n}\left( {0,1} \right),$ ${U_1},{U_2},...,{U_N}$ are subsets of ${B^n}\left( {0,R} \right) \times \left[ {0,R} \right],$ and $1 \le p,q <  + \infty ,$ $\Pi $ is a linear sub-space in $\mathbb{R}^{n+1}$ with dimension $m$, $1 \le m \le n + 1,$ $\pi $ is the orthogonal projection from $\mathbb{R}^{n+1}$ to  $\Pi $, then there exists a non-zero polynomial ${P_\Pi }$ defined on $\Pi $ of degree no more than ${C_m}{N^{1/m}},$ such that $P\left( z \right) = {P_\Pi }\left( {\pi \left( z \right)} \right),$ $z \in {\mathbb{R}^{n + 1}},$ satisfies
\begin{equation}\label{Eq41}
\left\| {{e^{it{\rm H}}}{f_j}{\psi _2}\left( {\frac{t}{R}} \right)} \right\|_{BL_{k,A}^p{L^q}\left( {{U_j} \cap \left\{ {P > 0} \right\}} \right)}^p = \left\| {{e^{it{\rm H}}}{f_j}{\psi _2}\left( {\frac{t}{R}} \right)} \right\|_{BL_{k,A}^p{L^q}\left( {{U_j} \cap \left\{ {P < 0} \right\}} \right)}^p,j = 1,2,...,N.
\end{equation}
\end{lemma}
\textbf{Proof:}  Let $V = \left\{ {P\left( z \right) = {P_\Pi }\left( {\pi \left( z \right)} \right):Deg{P_\Pi } \le D} \right\},$ note that $V$ is a vector space of dimension ${D^m},$ choose $D$ such that ${D^m} \sim N + 1,$ i.e., $D \le {C_m}{N^{1/m}},$ without less of generality, we may assume $DimV = N + 1$ and identify $V$ with $\mathbb{R}^{N+1}$. We define a function $G:{S^N} \subset V\backslash \left\{ 0 \right\} \to {\mathbb{R}^N}$ as
\[G\left( P \right) := \left\{ {{G_j}\left( P \right)} \right\}_{j = 1}^N,\]
where
\[{G_j}\left( P \right) := \left\| {{e^{it{\rm H}}}{f_j}{\psi _2}\left( {\frac{t}{R}} \right)} \right\|_{BL_{k,A}^p{L^q}\left( {{U_j} \cap \left\{ {P > 0} \right\}} \right)}^p - \left\| {{e^{it{\rm H}}}{f_j}{\psi _2}\left( {\frac{t}{R}} \right)} \right\|_{BL_{k,A}^p{L^q}\left( {{U_j} \cap \left\{ {P < 0} \right\}} \right)}^p.\]

It is obvious that $G\left( { - P} \right) =  - G\left( P \right).$ If the function $G$ is continuous, then Lemma \ref{Lemma4.1}  follows from the Borsuk - Ulam Theorem. So we only need to check the continuity of ${G_j}.$

Suppose ${P_l} \to P$ in $V\backslash \left\{ 0 \right\},$ note that
\[\left| {{G_j}\left( {{P_l}} \right) - {G_j}\left( P \right)} \right| \le 2\left\| {{e^{it{\rm H}}}{f_j}{\psi _2}\left( {\frac{t}{R}} \right)} \right\|_{BL_{k,A}^p{L^q}\left( {{U_j} \cap \left\{ {P{P_l} \le 0} \right\}} \right)}^p,\]
so we have
\[\mathop {\lim }\limits_{l \to  + \infty } \left\| {{e^{it{\rm H}}}{f_j}{\psi _2}\left( {\frac{t}{R}} \right)} \right\|_{BL_{k,A}^p{L^q}\left( {{U_j} \cap \left\{ {P{P_l} \le 0} \right\}} \right)}^p \le \left\| {{e^{it{\rm H}}}{f_j}{\psi _2}\left( {\frac{t}{R}} \right)} \right\|_{BL_{k,A}^p{L^q}\left( {{U_j} \cap {P^{ - 1}}\left( 0 \right)} \right)}^p = 0.\]
This implies that $G$ is continuous on $V\backslash \left\{ 0 \right\}.$

We use this Lemma to prove the following partitioning result:
\begin{theorem}\label{Theorem4.2}
 Suppose that $f$ is a function defined on $\mathbb{R}^{n}$ with $supp \widehat{f} \subset {B^n}\left( {0,1} \right),$ $U$ is a subset of ${B^n}\left( {0,R} \right) \times \left[ {0,R} \right],$ and $1 \le p,q <  + \infty ,$ $\Pi $ is a linear sub-space in $\mathbb{R}^{n+1}$ with dimension $m$, $1 \le m \le n + 1$, $\pi $ is the orthogonal projection from $\mathbb{R}^{n+1}$ to $\Pi $, then there exists a non-zero polynomial ${P_\Pi }$ defined on $\Pi $ of degree no more than $D$, and $P\left( z \right) = {P_\Pi }\left( {\pi \left( z \right)} \right)$
 such that $\Pi $ is a union of  ${ \sim _m}{D^{m}}$ disjoint open sets ${{\rm O}_{\Pi ,i}}$, ${\mathbb{R}^{n + 1}}\backslash Z\left( P \right)$ is a union of ${ \sim _m}{D^{m}}$ disjoint open sets ${{\rm O}_i} = {\pi ^{ - 1}}\left( {{{\rm O}_{\Pi ,i}}} \right),$ and for each $i$, we have
\begin{equation}\label{Eq42}
\left\| {{e^{it{\rm H}}}f{\psi _2}\left( {\frac{t}{R}} \right)} \right\|_{BL_{k,A}^p{L^q}\left( U \right)}^p \le {C_m}{D^m}\left\| {{e^{it{\rm H}}}f{\psi _2}\left( {\frac{t}{R}} \right)} \right\|_{BL_{k,A}^p{L^q}\left( {U \cap {{\rm O}_i}} \right)}^p.
\end{equation}
\end{theorem}
\textbf{Proof:} By Lemma \ref{Lemma4.1}, we obtain a polynomial ${Q_1}$ of degree $ \le C,$
\[{Q_1}\left( z \right) = {Q_{\Pi ,1}}\left( {\pi \left( z \right)} \right),\]
such that
\[\left\| {{e^{it{\rm H}}}f{\psi _2}\left( {\frac{t}{R}} \right)} \right\|_{BL_{k,A}^p{L^q}\left( {U \cap \left\{ {{Q_1} > 0} \right\}} \right)}^p = \left\| {{e^{it{\rm H}}}f{\psi _2}\left( {\frac{t}{R}} \right)} \right\|_{BL_{k,A}^p{L^q}\left( {U \cap \left\{ {{Q_1} < 0} \right\}} \right)}^p.\]
Next by Lemma \ref{Lemma4.1} again, we have a polynomial ${Q_2}$ of degree $ \le {C_m}{2^{1/m}},$ such that
\[\left\| {{e^{it{\rm H}}}f{\psi _2}\left( {\frac{t}{R}} \right)} \right\|_{BL_{k,A}^p{L^q}\left( {U \cap \left\{ {{Q_1} > 0} \right\} \cap \left\{ {{Q_2} > 0} \right\}} \right)}^p = \left\| {{e^{it{\rm H}}}f{\psi _2}\left( {\frac{t}{R}} \right)} \right\|_{BL_{k,A}^p{L^q}\left( {U \cap \left\{ {{Q_1} > 0} \right\} \cap \left\{ {{Q_2} < 0} \right\}} \right)}^p,\]
\[\left\| {{e^{it{\rm H}}}f{\psi _2}\left( {\frac{t}{R}} \right)} \right\|_{BL_{k,A}^p{L^q}\left( {U \cap \left\{ {{Q_1} < 0} \right\} \cap \left\{ {{Q_2} > 0} \right\}} \right)}^p = \left\| {{e^{it{\rm H}}}f{\psi _2}\left( {\frac{t}{R}} \right)} \right\|_{BL_{k,A}^p{L^q}\left( {U \cap \left\{ {{Q_1} < 0} \right\} \cap \left\{ {{Q_2} < 0} \right\}} \right)}^p.\]
Continuing inductively, we construct polynomials ${Q_{{1_{}}}},{Q_2},...,{Q_s},$
\[{Q_l}\left( z \right) = {Q_{\Pi ,l}}\left( {\pi \left( z \right)} \right),l = 1,2,...,s.\]
Set $P: = \prod\limits_{l = 1}^s {{Q_s}} ,$ where $\deg {Q_l} \le {C_m}{2^{\left( {l - 1} \right)/m}},$ therefore $\deg P\left( z \right) \le {C_m}{2^{s/m}},$ and the sign conditions of polynomials cut ${\mathbb{R}^{n + 1}}\backslash Z\left( P \right)$ into ${2^s}$ cells ${{\rm O}_i}$ such that
\[\left\| {{e^{it{\rm H}}}f{\psi _2}\left( {\frac{t}{R}} \right)} \right\|_{BL_{k,A}^p{L^q}\left( U \right)}^p \le {C_m}{2^s}\left\| {{e^{it{\rm H}}}f{\psi _2}\left( {\frac{t}{R}} \right)} \right\|_{BL_{k,A}^p{L^q}\left( {U \cap {{\rm O}_i}} \right)}^p.\]
Choose $s$ such that ${2^{s/m}} \in \left[ {D/2,D} \right],$ then we have $\deg P \le D$ and the number of cells ${{\rm O}_i}$ is ${C_m}{D^m}.$ It is obvious that $\Pi $ is divided by ${C_m}{D^m}$ cells ${{\rm O}_{\Pi ,i}}$ determined by the sign conditions of ${Q_{\Pi ,l}},l = 1,2,...,s.$ This completes the proof of Theorem \ref{Theorem4.2}.

Same as the analysis in \cite{DL}, by a slight modification in Theorem \ref{Theorem4.2}, we assume that all the varieties appear in our argument are transverse complete intersections. For any $1 \le m \le n,$ we say that a variety $Z\left( {{P_1},{P_2},...,{P_{n + 1 - m}}} \right) \subset {\mathbb{R}^n} \times \mathbb{R}$ is a transverse complete intersection if  for each $z \in Z\left( {{P_1},{P_2},...,{P_{n + 1 - m}}} \right)$, $\nabla {P_1}\left( z \right)\wedge \nabla {P_2}\left( z \right)\wedge...\wedge\nabla {P_{n + 1 - m}}\left( z \right)\neq0$.

  \section{Proof of Theorem \ref{Theorem3.2}}\label{Section 5}
P\textbf{roof of Theorem \ref{Theorem3.2}.} By Lemma \ref{Lemma3.3} (3), it is sufficient to prove Theorem \ref{Theorem3.2} for $p = 3.2.$ We assume that (\ref{Eq31}) holds for $A \le \frac{{\overline {A} }}{2}$ and ${R^{'}} \le \frac{R}{2},$ next prove it for $A = \overline {A} $ and ${R^{'}} = R.$

We say that we are in the algebraic case if there is a transverse complete intersection $Z\left( P \right)$ of dimension $2$, where $\deg P\left( z \right) \le D = D\left( \varepsilon  \right),$ so that
\begin{equation}\label{Eq51}
\left\| {{e^{it{\rm H}}}f{\psi _2}\left( {\frac{t}{R}} \right)} \right\|_{BL_{k,A}^p{L^q}\left( {B\left( {0,R} \right) \times \left[ {0,R} \right]} \right)}^{} \le C\left\| {{e^{it{\rm H}}}f{\psi _2}\left( {\frac{t}{R}} \right)} \right\|_{BL_{k,A}^p{L^q}(\left( {B\left( {0,R} \right) \times \left[ {0,R} \right]) \cap {N_{{R^{1/2 + \delta }}}}\left( {Z\left( P \right)} \right)} \right)}^{},
\end{equation}
here $N_{R^{1/2+\delta}}(Z(P))$ denotes the $R^{1/2+\delta}$ neighborhood of $Z(P)$. Otherwise we are in the cellular case.

\textbf{Cellular case.} We will use polynomial partitioning. By Theorem \ref{Theorem4.2}, there exists a non-zero polynomial $P\left( z \right) = \prod\limits_l {{Q_l}\left( z \right)} $ of degree at most $D$ such that $\left( {{\mathbb{R}^2} \times \mathbb{R}} \right)\backslash Z\left( P \right)$ is a union of $\sim {D^3}$ disjoint cells ${{\rm O}_i}$ such that for each $i$, we have
\[\left\| {{e^{it{\rm H}}}f{\psi _2}\left( {\frac{t}{R}} \right)} \right\|_{BL_{k,A}^p{L^q}\left( {B\left( {0,R} \right) \times \left[ {0,R} \right]} \right)}^p \le C{D^3}\left\| {{e^{it{\rm H}}}f{\psi _2}\left( {\frac{t}{R}} \right)} \right\|_{BL_{k,A}^p{L^q}\left( {\left( {B\left( {0,R} \right) \times \left[ {0,R} \right]} \right) \cap {{\rm O}_i}} \right)}^p.\]
Moreover, $Z\left( P \right)$ is a transverse complete intersection of dimension $2$.

Put
\[W := {N_{{R^{1/2 + \delta }}}}\left( {Z\left( P \right)} \right),{\rm O}_i^{'} := {\rm O}_i^{}\backslash W.\]
Since we are in the cellular case and $W \subset  \cup {}_l{N_{{R^{1/2 + \delta }}}}\left( {Z\left( {{Q_l}} \right)} \right),$ the contribution from $W$ is negligible. Hence for each $i$,
\begin{equation}\label{Eq52}
\left\| {{e^{it{\rm H}}}f{\psi _2}\left( {\frac{t}{R}} \right)} \right\|_{BL_{k,A}^p{L^q}\left( {B\left( {0,R} \right) \times \left[ {0,R} \right]} \right)}^p \le C{D^3}\left\| {{e^{it{\rm H}}}f{\psi _2}\left( {\frac{t}{R}} \right)} \right\|_{BL_{k,A}^p{L^q}\left( {\left( {B\left( {0,R} \right) \times \left[ {0,R} \right]} \right) \cap {\rm O}_i^{'}} \right)}^p.
\end{equation}

For each cell ${\rm O}_i^{'},$ we set
\[{{\rm T}_i} := \left\{ {\left( {\theta ,\upsilon } \right) \in {\rm T}:{T_{\theta ,\upsilon }} \cap {\rm O}_i^{'} \ne \emptyset} \right\}.\]
For the function $f$, we define
\[{f_i}: = \sum\limits_{\left( {\theta ,\upsilon } \right) \in {{\rm T}_i}} {{f_{\theta ,\upsilon }}} .\]
It follows that on ${\rm O}_i^{'},$
\[{e^{it{\rm H}}}f \sim {e^{it{\rm H}}}{f_i}.\]
By the fundamental theorem of Algebra, see \cite{DL}, for each $\left( {\theta ,\upsilon } \right) \in {\rm T},$ we have
\[Card\left\{ {i:\left( {\theta ,\upsilon } \right) \in {{\rm T}_i}} \right\} \le D + 1.\]
Hence
\[\sum\limits_i {\left\| {{f_i}} \right\|_{{L^2}}^2}  \le CD\left\| f \right\|_{{L^2}}^2,\]
by pigeonhole principle, there exists ${\rm O}_i^{'}$ such that
\begin{equation}\label{Eq53}
\left\| {{f_i}} \right\|_{{L^2}}^2 \le C{D^{ - 2}}\left\| f \right\|_{{L^2}}^2.
\end{equation}
So for such $i$, by (\ref{Eq52}), the induction on $R^{'}$ and (\ref{Eq53}), we have
\begin{align}
\left\| {{e^{it{\rm H}}}f{\psi _2}\left( {\frac{t}{R}} \right)} \right\|_{BL_{k,\overline{A}}^p{L^q}\left( {B\left( {0,R} \right) \times \left[ {0,R} \right]} \right)}^p &\le C{D^3}\left\| {{e^{it{\rm H}}}{f_i}{\psi _2}\left( {\frac{t}{R}} \right)} \right\|_{BL_{k,\overline{A}}^p{L^q}\left( {\left( {B\left( {0,R} \right) \times \left[ {0,R} \right]} \right) \cap {\rm O}_i^{'}} \right)}^p  \nonumber\\
  &\le C{D^3}\left\| {{e^{it{\rm H}}}{f_i}{\psi _2}\left( {\frac{t}{R}} \right)} \right\|_{BL_{k,\overline{ A} }^p{L^q}\left( {B\left( {0,R} \right) \times \left[ {0,R} \right]} \right)}^p  \nonumber\\
  &\le C{D^3}\sum\limits_{{B_{R/2}}\;{\mathop{\rm cov}} er\;B\left( {0,R} \right)} {\left\| {{e^{it{\rm H}}}{f_i}{\psi _2}\left( {\frac{t}{R}} \right)} \right\|_{BL_{k,\overline {A} }^p{L^q}\left( {{B_{R/2}} \times \left[ {0,R} \right]} \right)}^p}  \nonumber\\
  &\le C{D^{3 - p}}{\left( {C\left( {K,\varepsilon } \right){R^{\frac{1}{{qp}}}}{R^\varepsilon }{{\left\| f \right\|}_{{L^2}}}} \right)^p}, \nonumber
 \end{align}
choose $D$ sufficiently large such that $C{D^{3 - p}} \ll 1,$ this completes the induction.

\textbf{ Algebraic case.} We decompose $B\left( {0,R} \right) \times \left[ {0,R} \right]$ into balls ${B_j}$ of radius $\rho ,$ ${\rho ^{1/2 + {\delta _2}}} = {R^{1/2 + \delta }}.$ Choose $\delta  \ll {\delta _2},$ so that $\rho  \sim {R^{1 - O\left( {{\delta _2}} \right)}}.$ For each $j$ we define
\[{{\rm T}_j}: = \left\{ {\left( {\theta ,\upsilon } \right) \in {\rm T}:{T_{\theta ,\upsilon }} \cap {N_{{R^{1/2 + \delta }}}}\left( {Z\left( P \right)} \right) \cap {B_j} \ne \emptyset } \right\},\]
and
\[{f_j}: = \sum\limits_{\left( {\theta ,\upsilon } \right) \in {{\rm T}_j}} {{f_{\theta ,\upsilon }}} .\]
On each ${B_j} \cap {N_{{R^{1/2 + \delta }}}}\left( {Z\left( P \right)} \right)$, we have
\[{e^{it{\rm H}}}f \sim {e^{it{\rm H}}}{f_j}.\]
Therefore,
\[\left\| {{e^{it{\rm H}}}f{\psi _2}\left( {\frac{t}{R}} \right)} \right\|_{BL_{k,A}^p{L^q}\left( {B\left( {0,R} \right) \times [0,R]} \right)}^p \le \sum\limits_j {\left\| {{e^{it{\rm H}}}{f_j}{\psi _2}\left( {\frac{t}{R}} \right)} \right\|_{BL_{k,A}^p{L^q}\left( {{B_j} \cap {N_{{R^{1/2 + \delta }}}}\left( {Z\left( P \right)} \right)} \right)}^p} .\]

We further divide ${{\rm T}_j}$ into tubes that are tangential to $Z$ and tubes that are transverse to $Z$. We say that ${T_{\theta ,\upsilon }}$ is tangential to $Z$ in ${B_j}$ if the following two conditions hold:

\textbf{Distance condition:}
\[{T_{\theta ,\upsilon }} \cap 2{B_j} \subset {N_{{R^{1/2 + \delta }}}}\left( {Z\left( P \right)} \right) \cap 2{B_j} = {N_{{\rho ^{1/2 + {\delta _2}}}}}\left( {Z\left( P \right)} \right) \cap 2{B_j}.\]

\textbf{Angle condition:}
If $z \in Z \cap {N_{O\left( {{R^{1/2 + \delta }}} \right)}}\left( {{T_{\theta ,\upsilon }}} \right) \cap 2{B_j} = Z \cap {N_{O\left( {{\rho ^{1/2 + {\delta _2}}}} \right)}}\left( {{T_{\theta ,\upsilon }}} \right) \cap 2{B_j},$ then
\[Angle\left( {G\left( \theta  \right),{T_z}Z} \right) \le C{\rho ^{ - 1/2 + {\delta _2}}}.\]
The tangential wave packets are defined by
\[{{\rm T}_{j,{\rm{tang}}}} := \left\{ {\left( {\theta ,\upsilon } \right) \in {{\rm T}_j}:{T_{\theta ,\upsilon }} \text{ is tangent to } Z \text{ in } {B_j}} \right\},\]
and the transverse wave packets
\[{{\rm T}_{j,trans}}: = {{\rm T}_j}\backslash {{\rm T}_{j,{\rm{tang}}}}.\]
Set
\[{f_{j,{\rm{tang}}}}: = \sum\limits_{\left( {\theta ,\upsilon } \right) \in {{\rm T}_{j,{\rm{tang}}}}} {{f_{\theta ,\upsilon }}} , \hspace{0.2cm} {f_{j,{\rm{trans}}}}: = \sum\limits_{\left( {\theta ,\upsilon } \right) \in {{\rm T}_{j,{\rm{trans}}}}} {{f_{\theta ,\upsilon }}} ,\]
so
\[{f_j} = {f_{j,{\rm{tang}}}} + {f_{j,{\rm{trans}}}}.\]
Therefore, we have
\begin{align}
 \left\| {{e^{it{\rm H}}}f{\psi _2}\left( {\frac{t}{R}} \right)} \right\|_{BL_{k,A}^p{L^q}\left( {B\left( {0,R} \right) \times [0,R]} \right)}^p &\le \sum\limits_j {\left\| {{e^{it{\rm H}}}{f_j}{\psi _2}\left( {\frac{t}{R}} \right)} \right\|_{BL_{k,A}^p{L^q}\left( {{B_j}} \right)}^p}  \nonumber\\
  &\le \sum\limits_j {\left\| {{e^{it{\rm H}}}{f_{j,{\rm{tang}}}}{\psi _2}\left( {\frac{t}{R}} \right)} \right\|_{BL_{k,\frac{A}{2}}^p{L^q}\left( {{B_j}} \right)}^p} \nonumber\\
   \:\ &+ \sum\limits_j {\left\| {{e^{it{\rm H}}}{f_{j,{\rm{trans}}}}{\psi _2}\left( {\frac{t}{R}} \right)} \right\|_{BL_{k,\frac{A}{2}}^p{L^q}\left( {{B_j}} \right)}^p}.  \nonumber
 \end{align}
We will treat the tangential term and the transverse term respectively.

\textbf{ Algebraic transverse case.} In this case, the transverse term dominates, by induction on the radius $R^{'}$,
\begin{align}
 \left\| {{e^{it{\rm H}}}{f_{j,{\rm{trans}}}}{\psi _2}\left( {\frac{t}{R}} \right)} \right\|_{BL_{k,\frac{A}{2}}^p{L^q}\left( {{B_j}} \right)}^{} &\le \left\| {{e^{it{\rm H}}}{f_{j,{\rm{trans}}}}{\psi _2}\left( {\frac{t}{R}} \right)} \right\|_{BL_{k,\frac{A}{2}}^p{L^q}\left( {{B_\rho } \times \left[ {0,R} \right]} \right)}^{} \nonumber\\
  &\le C\left( {K,\varepsilon } \right){R^{\delta \left( {\log \overline A  - \log \overline {\frac{A}{2}} } \right)}}{R^{\frac{1}{{qp}}}}{\left( \rho  \right)^\varepsilon }{\left\| {{f_{j,{\rm{trans}}}}} \right\|_{{L^2}}} \nonumber\\
  &\le {R^{O\left( \delta  \right) - \varepsilon O\left( {{\delta _2}} \right)}}C\left( {K,\varepsilon } \right){R^{\frac{1}{{qp}}}}{R^\varepsilon }{\left\| {{f_{j,{\rm{trans}}}}} \right\|_{{L^2}}}, \nonumber
  \end{align}
where ${B_\rho }$ denotes the projection of $B_{j}$ on the $x$-plane. By \cite{DL} we have
\begin{equation}\label{Eq54}
\sum\limits_j {\left\| {{f_{j,{\rm{trans}}}}} \right\|_{{L^2}}^2}  \le C\left( D \right)\left\| f \right\|_{{L^2}}^2.
\end{equation}
Then
\begin{align}
 \sum\limits_j {\left\| {{e^{it{\rm H}}}{f_{j,{\rm{trans}}}}{\psi _2}\left( {\frac{t}{R}} \right)} \right\|_{BL_{k,\frac{A}{2}}^p{L^q}\left( {{B_j}} \right)}^p}  &\le {R^{O\left( \delta  \right) - \varepsilon O\left( {{\delta _2}} \right)}}{\left[ {C\left( {K,\varepsilon } \right){R^{\frac{1}{{qp}}}}{R^\varepsilon }} \right]^p}\sum\limits_j {\left\| {{f_{j,{\rm{trans}}}}} \right\|_{{L^2}}^p}  \nonumber\\
  &\le {R^{O\left( \delta  \right) - \varepsilon O\left( {{\delta _2}} \right)}}C\left( D \right){\left[ {C\left( {K,\varepsilon } \right){R^{\frac{1}{{qp}}}}{R^\varepsilon }\left\| f \right\|_{{L^2}}^{}} \right]^p}. \nonumber
\end{align}
The induction follows by choosing $\delta  \ll \varepsilon {\delta _2}$ and the fact that $R$ is sufficiently large.

\textbf{Algebraic tangential case.} in this case, the tangential term dominates, we need to do wave packets decomposition in ${B_j}$ at scale $\rho .$

\textbf{Wave packet decomposition in ${B_j}$.} Choose $\left( {\overline \theta  ,\overline {\nu}  } \right)$ as before where $\overline \theta$ is a ${\rho ^{ - 1/2}}$-cube in frequency space and $\overline \nu$ is a ${\rho ^{1/2}}$-cube in physical space. We can decompose $f$ as
\[f = \sum\limits_{\left( {\overline \theta  ,\overline \nu  } \right) \in {\rm T}} {{f_{\overline \theta  ,\overline \nu  }}}  = \sum\limits_{\left( {\overline \theta  ,\overline \nu  } \right) \in {\rm T}} {\left\langle {f,{\varphi _{\overline \theta  ,\overline \nu  }}} \right\rangle } {\varphi _{\overline \theta  ,\overline \nu  }},\]
where
\[\widehat{{\varphi _{\overline \theta  ,\overline \nu  }}}\left( \xi  \right) = {e^{ - ic\left( {\overline \nu  } \right) \cdot \xi }}{\widehat\varphi _{\overline \theta  }}\left( \xi  \right),\]
\[{\widehat\varphi _{\overline \theta  }}\left( {{\xi _1},{\xi _2}} \right) = \frac{1}{{{\rho ^{ - 1/2}}}}\prod\limits_{j = 1}^2 {\widehat\varphi } \left( {\frac{{{\xi _j} - c\left( {{\theta _j}} \right)}}{{{\rho ^{ - 1/2}}}}} \right).\]
Set $(x_{0},t_{0})$ as the center of $B_{j}$. In order to decompose wave packets in ${B_j}$, we need to modify the base such that
\begin{equation}\label{Eq55}
\hat{f} = \sum\limits_{\left( {\overline \theta  ,\overline \nu  } \right) \in {\rm T}} {\left\langle {\hat{f},{e^{ - i{x_0} \cdot \xi  + i\sqrt {{t_0}} \mu  \cdot \xi  - i{t_0}{{\left| \xi  \right|}^2}}}\widehat{{\varphi _{\overline \theta  ,\overline \nu  }}}\left( \xi  \right)} \right\rangle } {e^{ - i{x_0} \cdot \xi  + i\sqrt {{t_0}} \mu  \cdot \xi  - i{t_0}{{\left| \xi  \right|}^2}}}\widehat{{\varphi _{\overline \theta  ,\overline \nu  }}}\left( \xi  \right),
\end{equation}
so we set
\[\widehat{{{\tilde{\varphi}}_{\overline \theta  ,\overline \nu  }}} = {e^{ - i{x_0} \cdot \xi  + i\sqrt {{t_0}} \mu  \cdot \xi  - i{t_0}{{\left| \xi  \right|}^2}}}\widehat{{\varphi _{\overline \theta  ,\overline \nu  }}}\left( \xi  \right),\]
then
\begin{equation}\label{Eq56}
f = \sum\limits_{\left( {\overline \theta  ,\overline \nu  } \right) \in {\rm T}} {\left\langle {f,{{\tilde{\varphi}}_{\overline \theta  ,\overline \nu  }}} \right\rangle } {\tilde{\varphi}_{\overline \theta  ,\overline \nu  }}.
\end{equation}
Therefore,
\[{e^{it{\rm H}}}f = \sum\limits_{\left( {\overline \theta  ,\overline \nu  } \right) \in {\rm T}} {\left\langle {f,{{\tilde{\varphi}}_{\overline \theta  ,\overline \nu  }}} \right\rangle } {\tilde{\psi}_{\overline \theta  ,\overline \nu  }},\]
where
\[{\tilde{\psi}_{\overline \theta  ,\overline \nu  }} = {e^{it{\rm H}}}{\tilde{\varphi}_{\overline \theta  ,\overline \nu  }}.\]

As the previous analysis,  we restrict ${\tilde{\psi}_{\overline \theta  ,\overline \nu  }}$ in ${B_j}$, then we have
\[\left| {{{\tilde{\psi}}_{\overline \theta  ,\overline \nu  }}\left( {x,t} \right)} \right| \le {\rho ^{ - 1/2}}{\chi _{{T_{\overline \theta  ,\overline \nu  }}}}\left( {x,t} \right),\]
the tube ${T_{\overline \theta  ,\overline \nu  }}$ is defined by
\[{T_{\overline \theta  ,\overline \nu  }} := \left\{ {\left( {x,t} \right) \in {B_j}:\left| {x - {x_0} - c\left( {\overline \nu  } \right) + 2c\left( {\overline \theta  } \right)\left( {t - {t_0}} \right)} \right| \le {\rho ^{1/2 + \delta }},\left| {t - {t_0}} \right| \le \rho } \right\}.\]
For each $\left( {\theta ,\upsilon } \right) \in {{\rm T}_{j,{\rm{tang}}}},$ we consider the decomposition of ${f_{\theta ,\upsilon }},$
\[{f_{\theta ,\upsilon }} = \sum\limits_{\left( {\overline \theta  ,\overline \nu  } \right) \in {\rm T}} {\left\langle {{f_{\theta ,\upsilon }},{{\tilde{ \varphi}  }_{\overline \theta  ,\overline \nu  }}} \right\rangle } {\tilde{\varphi}_{\overline \theta  ,\overline \nu  }},\]
$\left( {\overline \theta  ,\overline \nu  } \right)$ which contribute to ${f_{\theta ,\upsilon }}$ satisfy
 \begin{equation}\label{Eq57}
\left| {c\left( \theta  \right) - c\left( {\overline \theta  } \right)} \right| \le 2{\rho ^{ - 1/2}},
\end{equation}
and
\begin{equation}\label{Eq58}
\left| {c\left( \nu  \right) - c\left( {\overline \nu  } \right) - {x_0}   - 2{t_0}c\left( \theta  \right)} \right| \le {R^{1/2 + \delta }}.
\end{equation}

From (\ref{Eq57}) we know that
 \begin{equation}\label{Eq59}
Angle\left( {G\left( \theta  \right),G\left( {\overline \theta  } \right)} \right) \le 2{\rho ^{ - 1/2}},
 \end{equation}
and (\ref{Eq58}) implies that if $\left( {x,t} \right) \in {T_{\overline \theta  ,\overline \nu  }},$ then
 \begin{equation}\label{Eq510}
\left| {x - c\left( \nu  \right) + 2c\left( \theta  \right)t} \right| \le C{R^{1/2 + \delta }},
 \end{equation}
i.e., ${T_{\overline \theta  ,\overline \nu  }} \subset {N_{{R^{1/2 + \delta }}}}\left( {{T_{\theta ,\upsilon }\cap B_{j}}} \right).$

We introduce the definition of ${\left( {{R^{'}}} \right)^{ - 1/2 + {\delta _m}}}$-tangent to $Z$ in $B$ with radius $R^{'}$.
Suppose that $Z = Z\left( {{P_1},...,{P_{3 - m}}} \right)$ is a transverse complete intersection in ${\mathbb{R}^2} \times \mathbb{R}.$ We say that ${T_{\theta ,\upsilon }}$(with scale ${R^{'}}$) is ${\left( {{R^{'}}} \right)^{ - 1/2 + {\delta _m}}}$-tangent to $Z$ in $B$ if the following two conditions hold:

 \textbf{(1) Distance condition:}
\[{T_{\theta ,\upsilon }} \subset {N_{{{\left( {{R^{'}}} \right)}^{1/2 + {\delta _m}}}}}\left( Z \right) \cap B.\]

\textbf{(2) Angle condition:} If $z \in Z \cap {N_{O\left( {{{\left( {{R^{'}}} \right)}^{1/2 + {\delta _m}}}} \right)}}\left( {{T_{\theta ,\upsilon }}} \right) \cap B,$ then
\[Angle\left( {G\left( \theta  \right),{T_z}Z} \right) \le C{\left( {{R^{'}}} \right)^{ - 1/2 + {\delta _m}}}.\]
Moreover, set
\[{{\rm T}_Z}: = \left\{ {\left( {\theta ,\upsilon } \right):{T_{\theta ,\upsilon }} \;is\; {\left( {{R^{'}}} \right)^{ - 1/2 + {\delta _m}}}{\rm{-tangent}}\;to\;Z\;in\;B} \right\},\]
we say that $f$ is concentrated in wave packets from ${{\rm T}_{Z}}$ in $B$ if
\[\sum\limits_{\left( {\theta ,\upsilon } \right) \notin {T_Z}} {\left\| {{f_{\theta ,\upsilon }}} \right\|_{{L^2}}^{}}  \le RapDec\left( {{R^{'}}} \right)\left\| f \right\|_{{L^2}}^{}.\]

We claim that new wave packets of ${f_{j,{\rm{tang}}}}$ are ${\rho ^{ - 1/2 + {\delta _2}}}$-tangent to $Z\left( P \right)$ in ${B_j}$ (note that we do not make a separate notation for convenience). In fact, if  $z \in Z \cap {N_{O\left( {{\rho ^{1/2 + {\delta _2}}}} \right)}}\left( {{T_{\overline \theta  ,\overline \nu  }}} \right) \cap {B_j},$ then $z \in Z \cap {N_{O\left( {{\rho ^{1/2 + {\delta _2}}}} \right)}}\left( {{T_{\theta ,\upsilon }}} \right) \cap {B_j},$ therefore
\[Angle\left( {G\left( {\overline \theta  } \right),{T_z}Z} \right) \le Angle\left( {G\left( {\overline \theta  } \right),G\left( \theta  \right)} \right) + Angle\left( {G\left( \theta  \right),{T_z}Z} \right) \le C\rho {^{-1/2 + {\delta _2}}}.\]
Also,
\[{T_{\overline \theta  ,\overline \nu  }} \subset {N_{{R^{1/2 + \delta }}}}\left( {{T_{\theta ,\upsilon }\cap B_{j}}} \right) \cap {B_j} = {N_{{\rho ^{1/2 + {\delta _2}}}}}\left( {{T_{\theta ,\upsilon }\cap B_{j} }} \right) \cap {B_j} \subset {N_{O\left( {{\rho ^{1/2 + {\delta _2}}}} \right)}}\left( {Z\left( P \right)} \right) \cap {B_j}.\]
Note that ${B_j} \subset {B_\rho } \times \left[ {0,R} \right],$ whenever Theorem \ref{Theorem6.1} below holds true, we have
\begin{align}
 &\left\| {{e^{it{\rm H}}}{f_{j,{\rm{tang}}}}{\psi _2}\left( {\frac{t}{R}} \right)} \right\|_{BL_{k,\frac{A}{2}}^p{L^q}\left( {{B_j}} \right)}^p \nonumber\\
 &\le {\left[ {{\rho ^{\left( {2 + 1/q} \right)\left( {1/p - 1/(4 + \delta) } \right)}}\left\| {{e^{it{\rm H}}}{f_{j,{\rm{tang}}}}{\psi _2}\left( {\frac{t}{R}} \right)} \right\|_{BL_{k,\frac{A}{2}}^{4 + \delta }{L^q}\left( {{B_j}} \right)}^{}} \right]^p} \nonumber\\
  &\le {\left[ {{\rho ^{\left( {2 + 1/q} \right)\left( {1/p - 1/(4 + \delta) } \right)}}\left\| {{e^{it{\rm H}}}{f_{j,{\rm{tang}}}}{\psi _2}\left( {\frac{t}{R}} \right)} \right\|_{BL_{k,\frac{A}{2}}^{4 + \delta }{L^q}\left( {{B_\rho } \times \left[ {0,R} \right]} \right)}^{}} \right]^p} \nonumber\\
  &\le {\left[ {{\rho ^{\left( {2 + 1/q} \right)\left( {1/p - 1/(4 + \delta) } \right)}}C\left( {K,D,\frac{\varepsilon }{2}} \right){R^{\delta \left( {\log \overline A  - \log \overline A /2} \right)}}{R^{\frac{1}{{qp}}}}{{\left( \rho  \right)}^{\frac{1}{{2\left( {4 + \delta } \right)}} - \frac{1}{4} + \frac{\varepsilon }{2}}}{{\left\| {{f_{j,{\rm{tang}}}}} \right\|}_{{L^2}}}} \right]^p} \nonumber
   \end{align}
  \begin{align}
  &\le {\left[ {{R^{O\left( \delta  \right) - \varepsilon /2}}C\left( {K,D,\frac{\varepsilon }{2}} \right){R^{\frac{1}{{qp}}}}{R^\varepsilon }{{\left\| {{f_{j,{\rm{tang}}}}} \right\|}_{{L^2}}}} \right]^p} \nonumber\\
  &\le {R^{O\left( \delta  \right) - p\varepsilon /2}}{\left[ {C\left( {K,\varepsilon } \right){R^{\frac{1}{{qp}}}}{R^\varepsilon }{{\left\| f \right\|}_{{L^2}}}} \right]^p}, \nonumber
 \end{align}
where we choose $C\left( {K,\varepsilon } \right) \ge C\left( {K,D,\frac{\varepsilon }{2}} \right),$ therefore,
\[\sum\limits_j {\left\| {{e^{it{\rm H}}}{f_{j,{\rm{tang}}}}{\psi _2}\left( {\frac{t}{R}} \right)} \right\|_{BL_{k,\frac{A}{2}}^p{L^q}\left( {{B_j}} \right)}^p}  \le {R^{O\left( {{\delta _2}} \right)}}{R^{O\left( \delta  \right) - p\varepsilon /2}}{\left[ {C\left( {K,\varepsilon } \right){R^{\frac{1}{{qp}}}}{R^\varepsilon }{{\left\| f \right\|}_{{L^2}}}} \right]^p},\]
the induction closes for the fact that ${\delta} \ll {\delta _2} \ll \varepsilon $ and $R$ is sufficiently large.

\begin{theorem}\label{Theorem6.1}
Suppose that $Z\left( P \right) \subset {\mathbb{R}^2} \times \mathbb{R}$ is a transverse complete intersection determined by some $P\left( z \right)$ with $degP\left( z \right) \le D_{Z}.$ For all $f$ with $ supp \hat{f} \subset B\left( {0,1} \right)$, and fixed $R \ge 1,$ if $B\left( {0,{R^{'}}} \right) \times \left[ {0,R} \right]$ contains a ball  (tube) $B$ of radius ${R^{'}}$ such that $f$ is concentrated in wave packets from ${T_Z}$ in $B$, here $1 \le {R^{'}} \le R,$ then  for any $\varepsilon  > 0$ and $p > 4,$ there exist positive constants $\overline A  = \overline A \left( \varepsilon  \right)$ and $C\left( {K,D_{Z},\varepsilon } \right)$ such that
\begin{equation}\label{Eq61}
\left\| {{e^{it{\rm H}}}f{\psi _2}\left( {\frac{t}{R}} \right)} \right\|_{BL_{k,A}^p{L^q}\left( {B\left( {0,{R^{'}}} \right) \times [0,R]} \right)}^{} \le C\left( {K,D_{Z},\varepsilon } \right){R^{\delta \left( {\log \overline A  - \log A} \right)}}{R^{\frac{1}{{qp}}}}{\left( {{R^{'}}} \right)^{\frac{1}{{2p}} - \frac{1}{4} + \varepsilon }}{\left\| f \right\|_{{L^2}}}
\end{equation}
holds for all $1 \le A \le \overline A .$
\end{theorem}

  \section{Proof of Theorem \ref{Theorem6.1}}\label{Section 6}
We will again use the induction on $R^{'}$ and $A$ to prove Theorem \ref{Theorem6.1}, the base of the induction is done as in Section 3. And we only consider the case $KM < {R^{1/2 - {\rm O}\left( {{\delta _1}} \right)}}$. We assume that the result holds for $A \le \frac{{\overline A }}{2}$ and ${R^{'}} \le \frac{R}{2},$ next prove it for  $A = \overline A $ and ${R^{'}} = R,$ this completes the induction.

Set $D = D\left( {\varepsilon, {D_Z}} \right)$, we say we are in algebraic case if there is transverse complete intersection ${Y} \subset {Z}$ of dimension $1$ defined using polynomials of degree no more than $D$, such that
\[\left\| {{e^{it{\rm H}}}f{\psi _2}\left( {\frac{t}{R}} \right)} \right\|_{BL_{k,A}^p{L^q}\left( {B\left( {0,R} \right) \times \left[ {0,R} \right]} \right)}^{} \le C\left\| {{e^{it{\rm H}}}f{\psi _2}\left( {\frac{t}{R}} \right)} \right\|_{BL_{k,A}^p{L^q}\left( ({B\left( {0,R} \right) \times \left[ {0,R} \right]) \cap {N_{{R^{1/2 + \delta_{2} }}}}\left( Y \right)} \right)}^{}.\]
Otherwise we are in the cellular case.

\textbf{Cellular case.} We first identify a significant piece $N_{1}$ of $(B\left( {0,R} \right) \times \left[ 0,R \right]) \cap {N_{{R^{1/2 + {\delta _2}}}}}\left( {Z\left( P \right)} \right)$,  where locally $Z\left( P \right)$ behaves like a $2$-plane $V$, such that
\begin{align}\label{Eq62}
 \left\| {{e^{it{\rm H}}}f{\psi _2}\left( {\frac{t}{R}} \right)} \right\|_{BL_{k,A}^p{L^q}\left( {B\left( {0,R} \right) \times \left[ {0,R} \right]} \right)}^{} &\le C\left\| {{e^{it{\rm H}}}f{\psi _2}\left( {\frac{t}{R}} \right)} \right\|_{BL_{k,A}^p{L^q}\left( ({B\left( {0,R} \right) \times \left[ {0,R} \right]) \cap {N_{{R^{1/2 + {\delta _2}}}}}\left( {Z\left( P \right)} \right)} \right)}^{} \nonumber\\
  &\le C\left\| {{e^{it{\rm H}}}f{\psi _2}\left( {\frac{t}{R}} \right)} \right\|_{BL_{k,A}^p{L^q}\left( {{N_1}} \right)}^{}.
 \end{align}
  By Theorem \ref{Theorem4.2}, there exists a polynomial $Q\left( z \right): = \prod\limits_{l = 1}^s {{Q_l}} $ with $\deg Q\left( z \right) \le D$, where polynomials ${Q_{{1_{}}}},{Q_2},...,{Q_s},$
\[{Q_l}\left( z \right) = {Q_{V,l}}\left( {\pi \left( z \right)} \right),l = 1,2,...,s,\]
$\pi $ is the orthogonal projection from ${\mathbb{R}^2} \times \mathbb{R}$ to $V,$ ${\mathbb{R}^2} \times R\backslash Z\left( Q \right)$ is divided into $\sim {D^2}$ cells ${{\rm O}_i}$ such that
 \begin{equation}\label{Eq63}
\left\| {{e^{it{\rm H}}}f{\psi _2}\left( {\frac{t}{R}} \right)} \right\|_{BL_{k,A}^p{L^q}\left( {{N_1}} \right)}^{p
} \le C{D^2}\left\| {{e^{it{\rm H}}}f{\psi _2}\left( {\frac{t}{R}} \right)} \right\|_{BL_{k,A}^p{L^q}\left( {{N_1} \cap {{\rm O}_i}} \right)}^p.
\end{equation}
 For each $l$, the variety ${Y_l} = Z\left( {P,{Q_l}} \right)$ is a transverse complete intersection of dimension $1$.  Define $W: = {N_{{R^{1/2 + \delta }}}}\left( {Z\left( Q \right)} \right),$ ${\rm O}_i^{'} := {\rm O}_i^{}\backslash W.$ By the analysis in \cite{G}, we have
 \[W \cap {N_1} \subset { \cup _l}{N_{O\left( {{R^{1/2 + \delta_{2} }}} \right)}}\left( {{Y_l}} \right),\]
 since we are in the cellular case, the contribution from $W$ is negligible. So we have
\begin{equation}\label{Eq64}
\left\| {{e^{it{\rm H}}}f{\psi _2}\left( {\frac{t}{R}} \right)} \right\|_{BL_{k,A}^p{L^q}\left( {{N_1}} \right)}^{p} \le C{D^2}\left\| {{e^{it{\rm H}}}f{\psi _2}\left( {\frac{t}{R}} \right)} \right\|_{BL_{k,A}^p{L^q}\left( {{N_1} \cap {\rm O}_i^{'}} \right)}^p.
\end{equation}
 Therefore,  from (\ref{Eq62})-(\ref{Eq64}) we actually obtain
 \begin{equation}\label{Eq65}
\left\| {{e^{it{\rm H}}}f{\psi _2}\left( {\frac{t}{R}} \right)} \right\|_{BL_{k,A}^p{L^q}\left( {B\left( {0,R} \right) \times [0,R]} \right)}^{p} \le C{D^2}\left\| {{e^{it{\rm H}}}f{\psi _2}\left( {\frac{t}{R}} \right)} \right\|_{BL_{k,A}^p{L^q}\left( ({B\left( {0,R} \right) \times [0,R]) \cap {\rm O}_i^{'}} \right)}^p.
\end{equation}

 For each cell ${\rm O}_i^{'},$ we set
\[{{\rm T}_i} := \left\{ {\left( {\theta ,\upsilon } \right) \in {\rm T}:{T_{\theta ,\upsilon }} \cap {\rm O}_i^{'} \ne \emptyset} \right\}.\]
 For the function $f$, we define
\[{f_i}: = \sum\limits_{\left( {\theta ,\upsilon } \right) \in {{\rm T}_i}} {{f_{\theta ,\upsilon }}} .\]
It follows that on ${\rm O}_i^{'},$
 \begin{equation}\label{Eq66}
{e^{it{\rm H}}}f \sim {e^{it{\rm H}}}{f_i}.
\end{equation}
 By the fundamental theorem of Algebra, for each $\left( {\theta ,\upsilon } \right) \in {\rm T},$
we have
\[Card\left\{ {i:£ý\left( {\theta ,\upsilon } \right) \in {{\rm T}_i}} \right\} \le D + 1.\]
Hence
\[\sum\limits_i {\left\| {{f_i}} \right\|_{{L^2}}^2}  \le CD\left\| f \right\|_{{L^2}}^2,\]
by pigeonhole principle, there exists ${\rm O}_i^{'}$ such that
\begin{equation}\label{Eq67}
\left\| {{f_i}} \right\|_{{L^2}}^2 \le C{D^{ - 1}}\left\| f \right\|_{{L^2}}^2.
\end{equation}
So by (\ref{Eq65}), (\ref{Eq66}), the induction on ${R^{'}},$ and (\ref{Eq67}), we have
\begin{align}
 \left\| {{e^{it{\rm H}}}f{\psi _2}\left( {\frac{t}{R}} \right)} \right\|_{BL_{k,\overline A }^p{L^q}\left( {B\left( {0,R} \right) \times [0,R]} \right)}^p &\le C{D^2}\left\| {{e^{it{\rm H}}}{f_i}{\psi _2}\left( {\frac{t}{R}} \right)} \right\|_{BL_{k,\overline A }^p{L^q}\left( {B\left( {0,R} \right) \times [0,R]} \right)}^p \nonumber\\
  &\le C{D^2}\sum\limits_{{B_{R/2}}\;{\mathop{\rm cov}} er\;B\left( {0,R} \right)} {\left\| {{e^{it{\rm H}}}{f_i}{\psi _2}\left( {\frac{t}{R}} \right)} \right\|_{BL_{k,\overline A }^p{L^q}\left( {{B_{R/2}} \times [0,R]} \right)}^p}  \nonumber\\
  &\le C{D^{2 - \frac{p}{2}}}{\left( {C\left( {K,D_{Z},\varepsilon } \right){R^{\frac{1}{{qp}}}}{{\left( {\frac{R}{2}} \right)}^{\frac{1}{{2p}} - \frac{1}{4} + \varepsilon }}{{\left\| f \right\|}_{{L^2}}}} \right)^p} \nonumber\\
  &\le C{D^{2 - \frac{p}{2}}}{\left( {C\left( {K,D_{Z},\varepsilon } \right){R^{\frac{1}{{qp}}}}{R^{\frac{1}{{2p}} - \frac{1}{4} + \varepsilon }}{{\left\| f \right\|}_{{L^2}}}} \right)^p}, \nonumber
 \end{align}
choose $D$ sufficiently large such that $C{D^{2 - \frac{p}{2}}} \ll 1,$ this completes the induction.

\textbf{Algebraic case. }In the algebraic case, there exists a transverse complete intersection $Y \subset Z\left( P \right)$ of dimension $1$, determined by polynomial with degree no more than  $D = D\left( {\varepsilon, {D_Z}} \right),$ so that
\[\left\| {{e^{it{\rm H}}}f{\psi _2}\left( {\frac{t}{R}} \right)} \right\|_{BL_{k,A}^p{L^q}\left( {B\left( {0,R} \right) \times [0,R]} \right)}^{} \le C\left\| {{e^{it{\rm H}}}f{\psi _2}\left( {\frac{t}{R}} \right)} \right\|_{BL_{k,A}^p{L^q}\left( ({B\left( {0,R} \right) \times  [0,R]) \cap {N_{{R^{1/2 + {\delta _2}}}}}\left( Y \right)} \right)}^{}.\]
We decompose $B\left( {0,R} \right) \times \left[ {0,R} \right]$ into balls ${B_j}$ of radius $\rho $, ${\rho ^{1/2 + {\delta _1}}} = {R^{1/2 + {\delta _2}}},$ ${\delta _2} \ll {\delta _1},$ in fact $\rho  \sim {R^{1 - O\left( {{\delta _1}} \right)}}.$ For each $j$, we define
\[{{\rm T}_j} := \left\{ {\left( {\theta ,\upsilon } \right) \in {\rm T}:{T_{\theta ,\upsilon }} \cap {N_{{R^{1/2 + {\delta _2}}}}}\left( {Y} \right) \cap {B_j} \ne \emptyset } \right\},\]
and
\[{f_j}: = \sum\limits_{\left( {\theta ,\upsilon } \right) \in {{\rm T}_j}} {{f_{\theta ,\upsilon }}} .\]
On each ${B_j} \cap {N_{{R^{1/2 + {\delta _2}}}}}(Y)$, we have
\[{e^{it{\rm H}}}f \sim {e^{it{\rm H}}}{f_j}.\]
Therefore,
\[\left\| {{e^{it{\rm H}}}f{\psi _2}\left( {\frac{t}{R}} \right)} \right\|_{BL_{k,A}^p{L^q}\left( {B\left( {0,R} \right) \times [0,R]} \right)}^p   \le \sum\limits_j {\left\| {{e^{it{\rm H}}}{f_j}{\psi _2}\left( {\frac{t}{R}} \right)} \right\|_{BL_{k,A}^p{L^q}\left( {{B_j}} \right)}^p} .\]
We further divide ${{\rm T}_j}$ into tubes that are tangential to $Y$ and tubes that are transverse to $Y$. We say that ${T_{\theta ,\upsilon }}$ is tangential to $Y$ in ${B_j}$ if the following two conditions hold:

\textbf{Distance condition:}
\begin{equation}\label{Eq68}
{T_{\theta ,\upsilon }} \cap 2{B_j} \subset {N_{{R^{1/2 + {\delta _2}}}}}\left( Y \right) \cap 2{B_j} = {N_{{\rho ^{1/2 + {\delta _1}}}}}\left( Y \right) \cap 2{B_j}.
\end{equation}

\textbf{Angle condition:}
If $z \in Y \cap {N_{O\left( {{R^{1/2 + {\delta _2}}}} \right)}}\left( {{T_{\theta ,\upsilon }}} \right) \cap 2{B_j} = Y \cap {N_{O\left( {{\rho ^{1/2 + {\delta _1}}}} \right)}}\left( {{T_{\theta ,\upsilon }}} \right) \cap 2{B_j},$ then
\begin{equation}\label{Eq69}
Angle\left( {G\left( \theta  \right),{T_z}Y} \right) \le C{\rho ^{ - 1/2 + {\delta _1}}}.
\end{equation}
The tangential wave packets is defined by
\[{{\rm T}_{j,{\rm{tang}}}} := \left\{ {\left( {\theta ,\upsilon } \right) \in {{\rm T}_j}:{T_{\theta ,\upsilon }} \text{ is tangent to } Y \text{ in } {B_j}} \right\},\]
and the transverse wave packets
\[{{\rm T}_{j,trans}}: = {{\rm T}_j}\backslash {{\rm T}_{j,{\rm{tang}}}}.\]
Set
\[{f_{j,{\rm{tang}}}}: = \sum\limits_{\left( {\theta ,\upsilon } \right) \in {{\rm T}_{j,{\rm{tang}}}}} {{f_{\theta ,\upsilon }}} , \hspace{0.2cm} {f_{j,{\rm{trans}}}}: = \sum\limits_{\left( {\theta ,\upsilon } \right) \in {{\rm T}_{j,{\rm{trans}}}}} {{f_{\theta ,\upsilon }}} ,\]
so
\[{f_j} = {f_{j,{\rm{tang}}}} + {f_{j,{\rm{trans}}}}.\]
Therefore, we have
\begin{align}
 \left\| {{e^{it{\rm H}}}f{\psi _2}\left( {\frac{t}{R}} \right)} \right\|_{BL_{k,A}^p{L^q}\left( {B\left( {0,R} \right) \times [0,R]} \right)}^p &\le \sum\limits_j {\left\| {{e^{it{\rm H}}}{f_j}{\psi _2}\left( {\frac{t}{R}} \right)} \right\|_{BL_{k,A}^p{L^q}\left( {{B_j}} \right)}^p}  \nonumber\\
  &\le \sum\limits_j {\left\| {{e^{it{\rm H}}}{f_{j,{\rm{tang}}}}{\psi _2}\left( {\frac{t}{R}} \right)} \right\|_{BL_{k,\frac{A}{2}}^p{L^q}\left( {{B_j}} \right)}^p} \nonumber\\
   &+ \sum\limits_j {\left\| {{e^{it{\rm H}}}{f_{j,{\rm{trans}}}}{\psi _2}\left( {\frac{t}{R}} \right)} \right\|_{BL_{k,\frac{A}{2}}^p{L^q}\left( {{B_j}} \right)}^p}.  \nonumber
 \end{align}
We will treat the tangential term and the transverse term respectively. Again, we need to use wave packets decomposition in ${B_j}.$

\textbf{ Algebraic tangential case.} In this case, the tangential term dominates. We claim that the new wave packets of ${f_{j,{\rm{tang}}}}$ are ${\rho ^{ - 1/2 + {\delta _1}}}$-tangent to $Y$ in ${B_j}$. In fact, by (\ref{Eq59}) and (\ref{Eq510}),
if $z \in Y\cap {N_{O\left( {{\rho ^{1/2 + {\delta _1}}}} \right)}}\left( {{T_{\overline \theta  ,\overline \nu  }}} \right) \cap {B_j},$ then $z \in Y \cap {N_{{O(R^{1/2 + {\delta _2}})}}}\left( {{T_{\theta ,\upsilon }}} \right) \cap {B_j},$ we have
\begin{equation}\label{Eq610}
Angle\left( {G\left( {\overline \theta  } \right),{T_z}Y} \right) \le Angle\left( {G\left( {\overline \theta  } \right),G\left( \theta  \right)} \right) + Angle\left( {G\left( \theta  \right),{T_z}Y} \right) \le C{\rho ^{ - 1/2 + {\delta _1}}}.
\end{equation}
Also,
\begin{equation}\label{Eq611}
{T_{\overline \theta  ,\overline \nu  }} \subset {N_{{R^{1/2 + {\delta _2}}}}}\left( {{T_{\theta ,\upsilon }\cap B_{j}}} \right) \cap {B_j} = {N_{{\rho ^{1/2 + {\delta _1}}}}}\left( {{T_{\theta ,\upsilon }\cap B_{j}}} \right) \cap {B_j} \subset {N_{O\left( {{\rho ^{1/2 + {\delta _1}}}} \right)}}\left( Y \right) \cap {B_j}.
\end{equation}
So, we can assume that ${f_{j,{\rm{tang}}}}$ is concentrated in wave packets from $T_{Y}$ in ${B_j}.$ Consider ${B_K} \times I_K^j$ such that
\[\left[ {{N_{O\left( {{\rho ^{1/2 + {\delta _1}}}} \right)}}\left( Y \right) \cap {B_j}} \right] \cap \left( {{B_K} \times I_K^j} \right) \ne \emptyset ,\]
there exists ${z_0} \in Y \cap {B_j} \cap {N_{O\left( {{\rho ^{1/2 + {\delta _1}}}} \right)}}\left( {{B_K} \times I_K^j} \right)$, for each ${T_{\overline \theta  ,\overline \nu  }}$ such that ${T_{\overline \theta  ,\overline \nu  }} \cap \left( {{B_K} \times I_K^j} \right) \ne \emptyset $, we have that ${z_0} \in Y \cap {B_j} \cap {N_{O\left( {{\rho ^{1/2 + {\delta _1}}}} \right)}}\left( {{T_{\overline \theta  ,\overline \nu  }}} \right)$, it holds
\[Angle\left( {G\left( {\overline \theta  } \right),{T_{{z_0}}}Y} \right) \le C{\rho ^{ - 1/2 + {\delta _1}}}.\]
Then for each $\tau $ with such a $\theta $ in it, it follows
\[Angle\left( {G\left( \tau  \right),{T_{{z_0}}}Y} \right) \le C{\rho ^{ - 1/2 + {\delta _1}}} \le {\left( {KM} \right)^{ - 1}},\]
such $\tau $ does not contribute to $\left\| {{e^{it{\rm H}}}{f_{j,{\rm{tang}}}}{\psi _2}\left( {\frac{t}{R}} \right)} \right\|_{BL_{k,\frac{A}{2}}^p{L^q}\left( {{B_j}} \right)}^p.$
 Since ${f_{j,{\rm{tang}}}}$ is concentrated in wave packets from $T_{Y}$ in ${B_j},$
\[\left\| {{e^{it{\rm H}}}{f_{j,{\rm{tang}}}}{\psi _2}\left( {\frac{t}{R}} \right)} \right\|_{BL_{k,\frac{A}{2}}^p{L^q}\left( {{B_j}} \right)}^p \le RapDec(\rho )\left\| f \right\|_{{L^2}}^p,\]
which can be negligible. So we only need to consider the transverse case.

\textbf{Algebraic transverse case.} In this case, the transverse term dominates. So we need to estimate
\[\sum\limits_j {\left\| {{e^{it{\rm H}}}{f_{j,{\rm{trans}}}}{\psi _2}\left( {\frac{t}{R}} \right)} \right\|_{BL_{k,\frac{A}{2}}^p{L^q}\left( {{B_j}} \right)}^p} .\]
Consider the new wave packets decomposition of ${f_{j,{\rm{trans}}}}$ in ${B_j},$ by (\ref{Eq59}) and (\ref{Eq510}), the new wave packets ${T_{\overline \theta  ,\overline \nu  }}$ satisfy
\begin{equation}\label{Eq612}
{T_{\overline \theta  ,\overline \nu  }} \subset {N_{{R^{1/2 + {\delta}}}}}\left( {{T_{\theta ,\upsilon }\cap B_{j}}} \right) \cap {B_j} \subset {N_{{R^{1/2 + {\delta _2}}}}}\left( Z \right) \cap {B_j}.
\end{equation}
And if $z \in Z \cap {N_{O\left( {{\rho^{1/2 + {\delta _2}}}} \right)}}\left( {{T_{\overline \theta  ,\overline \nu  }}} \right) \cap {B_j} \subset Z \cap {N_{O\left( {{R ^{1/2 + {\delta _2}}}} \right)}}\left( {{T_{\theta ,\upsilon }}} \right) \cap {B_j},$ then
\begin{equation}\label{Eq613}
Angle\left( {G\left( {\overline \theta  } \right),{T_z}Z} \right) \le Angle\left( {G\left( \theta  \right),{T_z}Z} \right) + Angle\left( {G\left( \theta  \right),G\left( {\overline \theta  } \right)} \right) \le C{\rho ^{ - 1/2 + {\delta _2}}}.
\end{equation}
So ${T_{\overline \theta  ,\overline \nu  }}$ is no longer $\rho^{-1/2+\delta_{2}}$-tangent to $Z$ in ${B_j}$ because the distance condition is not satisfied.

For each vector $b$ with $\left| b \right| \le {R^{1/2 + {\delta _2}}},$ define
\[{\overline {\rm T} _{Z + b}} := \left\{ {\left( {\overline \theta  ,\overline \upsilon  } \right):{T_{\overline \theta  ,\overline \upsilon  }} \text{ is } \rho^{-1/2+\delta_{2}} \text{-tangent to }Z + b \text{ in }{B_j}} \right\}.\]
By the angle condition, it turns out that each ${T_{\overline \theta  ,\overline \nu  }} \in {\overline {\rm T} _{Z + b}}$ for some $b$. We set
\[{f_{j,{\rm{trans,b}}}}: = \sum\limits_{\left( {\overline \theta  ,\overline \upsilon  } \right) \in {{\overline {\rm T} }_{Z + b}}} {{f_{\overline \theta  ,\overline \upsilon  }}} .\]
Then on ${B_j},$ it holds
\begin{equation}\label{Eq614}
\left| {{e^{it{\rm H}}}{f_{j,{\rm{trans,b}}}}} \right|{\psi _2}\left( {\frac{t}{R}} \right) \sim {\chi _{{N_{{\rho ^{1/2 + {\delta _2}}}}}\left( {Z + b} \right)}}\left( {x,t} \right)\left| {{e^{it{\rm H}}}{f_{j,{\rm{trans}}}}} \right|{\psi _2}\left( {\frac{t}{R}} \right).
\end{equation}

Next we choose a set of vectors $b \in {B_{R^{1/2 + {\delta _2}}}}.$ We cover ${N_{{R^{1/2 + {\delta _2}}}}}\left( Z \right) \cap {B_j}$ with disjoint balls of radius ${R^{1/2 + {\delta _2}}},$ and in each ball $B$ we note the value of ${N_{{\rho ^{1/2 + {\delta _2}}}}}\left( Z \right) \cap B.$ We will dyadically pigeonhole this volume.

For
\[{\mathcal{B}_s}: = \left\{ {B\left( {{x_0},{R^{1/2 + {\delta _2}}}} \right) \subset {N_{{R^{1/2 + {\delta _2}}}}}\left( Z \right) \cap {B_j}:B\left( {{x_0},{R^{1/2 + {\delta _2}}}} \right) \cap {N_{{\rho^{1/2 + {\delta _2}}}}}\left( Z \right) \sim {2^s}} \right\}.\]
We select a value of $s$ so that
\[\left\| {{e^{it{\rm H}}}{f_{j,{\rm{trans}}}}{\psi _2}\left( {\frac{t}{R}} \right)} \right\|_{BL_{k,\frac{A}{2}}^p{L^q}\left( {{B_j}} \right)}^p \le \left( {\log R} \right)\left\| {{e^{it{\rm H}}}{f_{j,{\rm{trans}}}}{\psi _2}\left( {\frac{t}{R}} \right)} \right\|_{BL_{k,\frac{A}{2}}^p{L^q}\left( {{ \cup _{B \in {{\mathcal{B}}_s}}B}} \right)}^p.\]
Therefore, we only consider $\left( {\theta ,\nu } \right)$ such that ${T_{\theta ,\upsilon }}$ meets at least one of the balls in $\mathcal{B}_s$. We choose a random set of $\left| {{B_{{R^{1/2 + {\delta _2}}}}}} \right|/{2^s}$ vectors $b \in {B_{{R^{1/2 + {\delta _2}}}}}$. For a typical ball $B\left( {{x_0},{R^{1/2 + {\delta _2}}}} \right) \in {\mathcal{B}_s},$ the union ${ \cup _b}{N_{{\rho ^{1/2 + {\delta _2}}}}}\left( {Z + b} \right) \cap {B_j}$ covers a definite fraction of the ball with high probability.  It follows
\begin{align}\label{Eq615}
&\left\| {{e^{it{\rm H}}}{f_{j,{\rm{trans}}}}{\psi _2}\left( {\frac{t}{R}} \right)} \right\|_{BL_{k,\frac{A}{2}}^p{L^q}\left( {{B_j}} \right)}^p \nonumber\\
&\le \left( {\log R} \right)\sum\limits_b {\left\| {{e^{it{\rm H}}}{f_{j,{\rm{trans,b}}}}{\psi _2}\left( {\frac{t}{R}} \right)} \right\|_{BL_{k,\frac{A}{2}}^p{L^q}\left( {{N_{{\rho ^{1/2 + {\delta _2}}}}}\left( {Z + b} \right) \cap {B_j}} \right)}^p}.
 \end{align}

 By the induction on $R^{'}$, we have
\begin{align}
 &\left\| {{e^{it{\rm H}}}{f_{j,{\rm{trans,b}}}}{\psi _2}\left( {\frac{t}{R}} \right)} \right\|_{BL_{k,\frac{A}{2}}^p{L^q}\left( {{N_{{\rho ^{1/2 + {\delta _2}}}}}\left( {Z + b} \right) \cap {B_j}} \right)}^p \nonumber\\
 &\le \left\| {{e^{it{\rm H}}}{f_{j,{\rm{trans,b}}}}{\psi _2}\left( {\frac{t}{R}} \right)} \right\|_{BL_{k,\frac{A}{2}}^p{L^q}\left( {{B_j}} \right)}^p \nonumber\\
  &\le \left\| {{e^{it{\rm H}}}{f_{j,{\rm{trans,b}}}}{\psi _2}\left( {\frac{t}{R}} \right)} \right\|_{BL_{k,\frac{A}{2}}^p{L^q}\left( {{B_\rho } \times \left[ {0,R} \right]} \right)}^p \nonumber\\
  &\le {\left[ {C\left( {K,{D_{Z}},\varepsilon } \right){R^{\delta \left( {\log \overline A  - \log \overline {\frac{A}{2}} } \right)}}{R^{\frac{1}{{qp}}}}{{\left( \rho  \right)}^{\frac{1}{{2p}} - \frac{1}{4} + \varepsilon }}{{\left\| {{f_{j,{\rm{trans,b}}}}} \right\|}_{{L^2}}}} \right]^p} \nonumber\\
  &\le {\left[ {C\left( {K,{D_Z},\varepsilon } \right){R^{O\left( \delta  \right)}}{R^{\frac{1}{{qp}}}}{{\left( \rho  \right)}^{\frac{1}{{2p}} - \frac{1}{4} + \varepsilon }}{{\left\| {{f_{j,{\rm{trans,b}}}}} \right\|}_{{L^2}}}} \right]^p}, \nonumber
  \end{align}
therefore, if
\begin{equation}\label{Eq621}
\sum\limits_j {\sum\limits_b {\left\| {{f_{j,{\rm{trans,b}}}}} \right\|_{{L^2}}^2} }  \sim \sum\limits_j {\left\| {{f_{j,{\rm{trans}}}}} \right\|_{{L^2}}^2}  \le D\left\| f \right\|_{{L^2}}^2,
\end{equation}
\begin{equation}\label{Eq620}
\mathop {\max }\limits_b \left\| {{f_{j,{\rm{trans,b}}}}} \right\|_{{L^2}}^2 \le {R^{O\left( {{\delta _2}} \right)}}{\left( {\frac{R}{\rho }} \right)^{ - 1/2}}\left\| {{f_{j,{\rm{trans}}}}} \right\|_{{L^2}}^2,
\end{equation}
then we have
\begin{align}
 &\sum\limits_j {\sum\limits_b {\left\| {{e^{it{\rm H}}}{f_{j,{\rm{trans}}}}{\psi _2}\left( {\frac{t}{R}} \right)} \right\|_{BL_{k,\frac{A}{2}}^p{L^q}\left( {{B_j}} \right)}^p} }  \nonumber\\
 &\le {\left[ {C\left( {K,{D_Z},\varepsilon } \right){R^{O\left( \delta  \right)}}{R^{\frac{1}{{qp}}}}{{\left( \rho  \right)}^{\frac{1}{{2p}} - \frac{1}{4} + \varepsilon }}} \right]^p}\sum\limits_j {\sum\limits_b {\left\| {{f_{j,{\rm{trans,b}}}}} \right\|_{{L^2}}^p} }  \nonumber\\
  &\le {\left[ {C\left( {K,{D_Z},\varepsilon } \right){R^{O\left( \delta  \right)}}{R^{\frac{1}{{qp}}}}{{\left( \rho  \right)}^{\frac{1}{{2p}} - \frac{1}{4} + \varepsilon }}} \right]^p}\,\sum\limits_j {\sum\limits_b {\left\| {{f_{j,{\rm{trans,b}}}}} \right\|_{{L^2}}^2\mathop {\mathop {\max }\limits_b \left\| {{f_{j,{\rm{trans,b}}}}} \right\|_{{L^2}}^{p - 2}}\limits_{} } }  \nonumber\\
  &\le {\left[ {C\left( {K,{D_Z},\varepsilon } \right){R^{O\left( \delta  \right)}}{R^{\frac{1}{{qp}}}}{{\left( \rho  \right)}^{\frac{1}{{2p}} - \frac{1}{4} + \varepsilon }}} \right]^p}\sum\limits_j {\sum\limits_b {\left\| {{f_{j,{\rm{trans,b}}}}} \right\|_{{L^2}}^2\mathop {\mathop {\max }\limits_b \left\| {{f_{j,{\rm{trans,b}}}}} \right\|_{{L^2}}^{p - 2}}\limits_{} } }  \nonumber\\
  &\le {\left[ {C\left( {K,{D_Z},\varepsilon } \right){R^{O\left( \delta  \right)}}{R^{\frac{1}{{qp}}}}{{\left( \rho  \right)}^{\frac{1}{{2p}} - \frac{1}{4} + \varepsilon }}} \right]^p}{R^{O\left( {{\delta _2}} \right)}}{\left( {\frac{R}{\rho }} \right)^{ - (p/2 - 1)}}\left\| f \right\|_{{L^2}}^p \nonumber\\
  &= {\left[ {C\left( {K,{D_Z},\varepsilon } \right){R^{O\left( \delta  \right)}}{R^{\frac{1}{{qp}}}}{R^{\frac{1}{{2p}} - \frac{1}{4} + \varepsilon  - O\left( {{\delta _1}} \right)\left( {\frac{1}{{2p}} - \frac{1}{4} + \varepsilon } \right)}}} \right]^p} \times {R^{O\left( {{\delta _2}} \right)}}{\left( {\frac{R}{\rho }} \right)^{ - (p/2 - 1)}}\left\| f \right\|_{{L^2}}^p  \nonumber\\
  &= {R^{O\left( \delta  \right)}}{R^{O\left( {{\delta _2}} \right)}}{R^{ - O\left( {{\delta _1}} \right)(p/2 - 1) - O\left( {{\delta _1}} \right)\left( {\frac{1}{2} - \frac{p}{4} + \varepsilon } \right)}}\;{\left[ {C\left( {K,{D_Z},\varepsilon } \right){R^{\frac{1}{{qp}}}}{R^{\frac{1}{{2p}} - \frac{1}{4} + \varepsilon }}\left\| f \right\|_{{L^2}}^{}} \right]^p} \nonumber\\
  &\le {R^{O\left( {{\delta _2}} \right)}}{R^{ - O\left( {{\delta _1}} \right)\varepsilon }}{\left[ {C\left( {K,{D_Z},\varepsilon } \right){R^{\frac{1}{{qp}}}}{R^{\frac{1}{{2p}} - \frac{1}{4} + \varepsilon }}\left\| f \right\|_{{L^2}}^{}} \right]^p}, \nonumber
\end{align}
so the induction closes by choosing ${\delta _2} \ll \varepsilon {\delta _1}$ and the fact that $R$ is sufficiently large. This completes the proof of Theorem \ref{Theorem6.1}.

Next we will prove (\ref{Eq621}) and (\ref{Eq620}). For each $\left( {\theta ,\nu } \right) \in {{\rm T}_{j,trans}},$ if ${T_{\theta ,\nu }}$ contributes, then ${T_{\theta ,\nu }}$ intersects some $B\left( {{x_0},{R^{1/2 + {\delta _2}}}} \right)$ in ${\mathcal{B}_s}.$ We have
\begin{equation}\label{Eq615}
\left\| {{f_{\theta ,\nu }}} \right\|_{{L^2}}^2 \sim {R^{ - 1/2 - {\delta _2}}}\left\| {{e^{it{\rm H}}}{f_{\theta ,\nu }}{\psi _2}\left( {\frac{t}{R}} \right)} \right\|_{{L^2}\left( {B\left( {{x_0},{R^{1/2 + {\delta _2}}}} \right)} \right)}^2,
\end{equation}
provided Theorem \ref{Theorem7.2} below holds true. Set
\begin{equation}\label{Eq616}
{f_{\theta ,\nu ,b}}: = \sum\limits_{\left( {\overline \theta  ,\overline \upsilon  } \right) \in {{\overline {\rm T} }_{Z + b}} \cap {{\left( {\theta ,\nu } \right)}^ \sim }} {{f_{\overline \theta  ,\overline \upsilon  }}} ,
\end{equation}
here $(\theta ,\nu)^{\sim}$  denotes the wave packets decomposition of $f_{\theta ,\nu}$ in $B_{j}$, it follows
\begin{equation}\label{Eq617}
\left\| {{e^{it{\rm H}}}{f_{\theta ,\nu ,b}}{\psi _2}\left( {\frac{t}{R}} \right)} \right\|_{{L^2}\left( {B\left( {{x_0},{R^{1/2 + {\delta _2}}}} \right)} \right)}^2 \sim \left\| {{e^{it{\rm H}}}{f_{\theta ,\nu }}{\psi _2}\left( {\frac{t}{R}} \right)} \right\|_{{L^2}\left( {B\left( {{x_0},{R^{1/2 + {\delta _2}}}} \right) \cap {N_{{\rho ^{1/2 + {\delta _2}}}}}\left( {Z + b} \right)} \right)}^2,
\end{equation}
using Theorem \ref{Theorem7.2} again,
\begin{equation}\label{Eq622}
\left\| {{f_{\theta ,\nu ,b}}} \right\|_{{L^2}}^2 \sim {R^{ - 1/2 - {\delta _2}}}\left\| {{e^{it{\rm H}}}{f_{\theta ,\nu ,b}}{\psi _2}\left( {\frac{t}{R}} \right)} \right\|_{{L^2}\left( {B\left( {{x_0},{R^{1/2 + {\delta _2}}}} \right)} \right)}^2.
\end{equation}
Notice that the sets ${N_{{\rho ^{1/2 + {\delta _2}}}}}\left( {Z + b} \right)$ are essentially disjoint, hence (\ref{Eq615}), (\ref{Eq616}) and (\ref{Eq622}) imply
\begin{align}
 \sum\limits_b {\left\| {{f_{\theta ,\nu ,b}}} \right\|_{{L^2}}^2}  &\sim {R^{ - 1/2 - {\delta _2}}}\sum\limits_b {\left\| {{e^{it{\rm H}}}{f_{\theta ,\nu }}{\psi _2}\left( {\frac{t}{R}} \right)} \right\|_{{L^2}\left( {B\left( {{x_0},{R^{1/2 + {\delta _2}}}} \right) \cap {N_{{\rho ^{1/2 + {\delta _2}}}}}\left( {Z + b} \right)} \right)}^2}  \nonumber\\
  &\le {R^{ - 1/2 - {\delta _2}}}\left\| {{e^{it{\rm H}}}{f_{\theta ,\nu }}{\psi _2}\left( {\frac{t}{R}} \right)} \right\|_{{L^2}\left( {B\left( {{x_0},{R^{1/2 + {\delta _2}}}} \right)} \right)}^2  \nonumber\\
  &\sim \left\| {{f_{\theta ,\nu }}} \right\|_{{L^2}}^2. \nonumber
 \end{align}
Therefore
\begin{equation}\label{Eq619}
\sum\limits_b {\left\| {{f_{j,{\rm{trans,b}}}}} \right\|_{{L^2}}^2}  = \sum\limits_{\left( {\theta ,\nu } \right)\in {{\rm T}_{j,trans}}} \sum\limits_b {\left\| {{f_{\theta ,\nu ,b}}} \right\|_{{L^2}}^2}  \le \sum\limits_{\left( {\theta ,\nu } \right) \in {{\rm T}_{j,trans}}} {\left\| {{f_{\theta ,\nu }}} \right\|_{{L^2}}^2}  = \left\| {{f_{j,{\rm{trans}}}}} \right\|_{{L^2}}^2.
\end{equation}
Then by (\ref{Eq54}) and (\ref{Eq619}), (\ref{Eq621}) holds.

If Theorem \ref{Theorem7.1} below holds true, then for each $b$,
\begin{align}
 \left\| {{f_{j,{\rm{trans,b}}}}} \right\|_{{L^2}}^2 &\le \sum\limits_{\left( {\theta ,\nu } \right) \in {{\rm T}_{j,trans}}} {} \left\| {{f_{\theta ,\nu ,b}}} \right\|_{{L^2}}^2 \sim {R^{ - 1/2 - {\delta _2}}}\sum\limits_{\left( {\theta ,\nu } \right) \in {{\rm T}_{j,trans}}} {\left\| {{e^{it{\rm H}}}{f_{\theta ,\nu ,b}}{\psi _2}\left( {\frac{t}{R}} \right)} \right\|_{{L^2}\left( {B\left( {{x_0},{R^{1/2 + {\delta _2}}}} \right)} \right)}^2}  \nonumber\\
  &\sim {R^{ - 1/2 - {\delta _2}}}\sum\limits_{\left( {\theta ,\nu } \right) \in {{\rm T}_{j,trans}}} {\left\| {{e^{it{\rm H}}}{f_{\theta ,\nu }}{\psi _2}\left( {\frac{t}{R}} \right)} \right\|_{{L^2}\left( {B\left( {{x_0},{R^{1/2 + {\delta _2}}}} \right) \cap {N_{{\rho ^{1/2 + {\delta _2}}}}}\left( {Z + b} \right)} \right)}^2,}  \nonumber
  \end{align}
 and
\begin{align}
&\left\| {{e^{it{\rm H}}}{f_{\theta ,\nu }}{\psi _2}\left( {\frac{t}{R}} \right)} \right\|_{{L^2}\left( {B\left( {{x_0},{R^{1/2 + {\delta _2}}}} \right) \cap {N_{{\rho ^{1/2 + {\delta _2}}}}}\left( {Z + b} \right)} \right)}^2 \nonumber\\
&\le C{R^{O\left( {{\delta _2}} \right)}}{\left( {\frac{{{R^{1/2}}}}{{{\rho ^{1/2}}}}} \right)^{ - 1}}\left\| {{e^{it{\rm H}}}{f_{\theta ,\nu }}{\psi _2}\left( {\frac{t}{R}} \right)} \right\|_{{L^2}\left( {B\left( {{x_0},2{R^{1/2 + {\delta _2}}}} \right)} \right)}^2,\nonumber
\end{align}
therefore
\begin{align}
 \left\| {{f_{j,{\rm{trans,b}}}}} \right\|_{{L^2}}^2 &\le C{R^{O\left( {{\delta _2}} \right)}}{\left( {\frac{{{R^{1/2}}}}{{{\rho ^{1/2}}}}} \right)^{ - 1}}{R^{ - 1/2 - {\delta _2}}}\sum\limits_{\left( {\theta ,\nu } \right)\in {{\rm T}_{j,trans}}} {\left\| {{e^{it{\rm H}}}{f_{\theta ,\nu }}{\psi _2}\left( {\frac{t}{R}} \right)} \right\|_{{L^2}\left( {B\left( {{x_0},2{R^{1/2 + {\delta _2}}}} \right)} \right)}^2}  \nonumber\\
  &\le C{R^{O\left( {{\delta _2}} \right)}}{\left( {\frac{{{R^{1/2}}}}{{{\rho ^{1/2}}}}} \right)^{ - 1}}\sum\limits_{\left( {\theta ,\nu } \right)\in {{\rm T}_{j,trans}}} {\left\| {{f_{\theta ,\nu }}} \right\|_{{L^2}}^2}  \nonumber\\
  &\le C{R^{O\left( {{\delta _2}} \right)}}{\left( {\frac{{{R^{1/2}}}}{{{\rho ^{1/2}}}}} \right)^{ - 1}}\left\| {{f_{j,{\rm{trans}}}}} \right\|_{{L^2}}^2, \nonumber
  \end{align}
in the second inequality above we used Theorem \ref{Theorem7.2} again, and (\ref{Eq620}) is obtained.

  \section{Transverse Equidistribution estimate}\label{Section 6}
The following Theorem \ref{Theorem7.1} is a generalization of Lemma 6.2 in \cite{G}, which is needed in the proof of Theorem \ref{Theorem6.1}. In order to prove Theorem \ref{Theorem7.1}, we need a version of the Heisenberg uncertainly principle in \cite{G}:

\begin{lemma}\label{Lemma7.1}
(\cite{G}) Suppose that $G:{\mathbb{R}^n} \to \mathbb{C}$ is a function, and that $\widehat{G}$ is supported in a ball $B\left( {{\xi _0},r} \right)$, then for any ball ${B_\rho }$ with $\rho  \le r^{-1},$ we have the inequality
\begin{equation}\label{Eq71}
\int_{{B_\rho }} {{{\left| G \right|}^2}}  \le C\frac{{\left| {{B_\rho }} \right|}}{{\left| {{B_{{r^{ - 1}}}}} \right|}}\int_{{B_{{r^{ - 1}}}}} {{{\left| G \right|}^2}}.
\end{equation}
\end{lemma}
\begin{theorem}\label{Theorem7.1}
 Suppose that $f$ is concentrated in wave packets from ${{\rm T}_Z},$ $Z = Z\left( P \right)$ is a  transverse complete intersection of dimension $2$, $B$ is a ball of radius ${R^{1/2 + {\delta _2}}}$ contained in $B\left( {0,R} \right) \times \left[ {0,R} \right],$ ${{\rm T}_{Z,B}} := \left\{ {\left( {\theta ,\nu } \right) \in {{\rm T}_Z}:{T_{\theta ,\nu }} \cap B \ne \emptyset } \right\}$, if
\[f = \sum\limits_{\left( {\theta ,\nu } \right) \in {{\rm T}_{Z,B}}} {{f_{\theta ,\nu }}} ,\]
then
\begin{equation}\label{Eq72}
\left\| {{e^{it{\rm H}}}f{\psi _2}\left( {\frac{t}{R}} \right)} \right\|_{{L^2}\left( {B \cap {N_{{\rho ^{1/2 + {\delta _2}}}}}\left( Z \right)} \right)}^2 \le C{R^{O\left( {{\delta _2}} \right)}}{\left( {\frac{{{R^{1/2}}}}{{{\rho ^{1/2}}}}} \right)^{ - 1}}\left\| {{e^{it{\rm H}}}f{\psi _2}\left( {\frac{t}{R}} \right)} \right\|_{{L^2}\left( {2B} \right)}^2.
\end{equation}
\end{theorem}
\textbf{Proof:} If $B \cap {N_{{R^{1/2 + {\delta _2}}}}}\left( Z \right) = \emptyset ,$ then ${{\rm T}_{Z,B}} = \emptyset ,$ and there is nothing to prove. So we can assume that $B \cap {N_{{R^{1/2 + {\delta _2}}}}}\left( Z \right) \ne \emptyset ,$ then there exists a point ${z_0} \in Z$ such that ${z_0} \in Z \cap {N_{{R^{1/2 + {\delta _2}}}}}\left( B \right)$, then for each wave packet $\left( {\theta ,\nu } \right) \in {{\rm T}_{Z,B}},$ we have \[{z_0} \in Z \cap {N_{{R^{1/2 + {\delta _2}}}}}\left( {{T_{\theta ,\nu }}} \right).\]
By the definition of ${{\rm T}_Z},$ we have
\begin{equation}\label{Eq73}
Angle\left( {G\left( \theta  \right),{T_{{z_0}}}Z} \right) \le {R^{ - 1/2 + {\delta _2}}}.
\end{equation}

We can assume ${T_{{z_0}}}Z$ is given by
\begin{equation}\label{Eq74}
{a_1}{x_1} + {a_2}{x_2} + bt = 0,\;a_1^2 + a_2^2 + {b^2} = 1,\;\left| {\left( {{a_1},{a_2}} \right)} \right| \ge 1,
\end{equation}
(\ref{Eq73}) and (\ref{Eq74}) imply
\[\left| { - 2c\left( \theta  \right) \cdot a + b} \right| \le C{R^{ - 1/2 + {\delta _2}}},\]
this restricts all $\theta $ to a strip of width ${R^{ - 1/2 + {\delta _2}}}$ paralleled to $\left( {{a_2}, - {a_1}} \right),$ we denote it by $S$.
The Fourier transform of ${e^{it{\rm H}}}{f_{\theta ,\nu }}{\psi _2}\left( {\frac{t}{R}} \right)$ is supported in
\[\left\{ {\left( {{\xi _1},{\xi _2},{\xi _3}} \right):\left( {{\xi _1},{\xi _2}} \right) \in \theta ,\left| {{\xi _3} - \xi _1^2 - \xi _2^2} \right| \le {R^{ - \varepsilon /2}}} \right\},\]
therefore, the Fourier transform of  ${e^{it{\rm H}}}f{\psi _2}\left( {\frac{t}{R}} \right)$ is supported in
\[\left\{ {\left( {{\xi _1},{\xi _2},{\xi _3}} \right):\left( {{\xi _1},{\xi _2}} \right) \in S,\left| {{\xi _3} - \xi _1^2 - \xi _2^2} \right| \le {R^{ - \varepsilon /2}}} \right\}.\]
Suppose that $\Pi $ is a $1$-dimension linear sub-space of ${\mathbb{R}^3}$ parallel to $\left( {{a_1}, {a_2},0} \right),$ then the projection of the Fourier transform of ${e^{it{\rm H}}}f{\psi _2}\left( {\frac{t}{R}} \right)$ on $\Pi $ is supported in a ball of radius ${R^{ - 1/2 + {\delta _2}}}$. If  we view ${e^{it{\rm H}}}f{\psi _2}\left( {\frac{t}{R}} \right)$ as a function defined on $\Pi $, then for each $x \in B \cap \Pi,$ Lemma \ref{Lemma7.1} implies
\begin{align}\label{Eq75}
 &\int_{\Pi  \cap B\left( {x,{\rho ^{1/2 + {\delta _2}}}} \right)} {{{\left| {{e^{it{\rm H}}}f{\psi _2}\left( {\frac{t}{R}} \right)} \right|}^2}}  \nonumber\\
 &\le \sum\limits_{B{}_{{\rho ^{1/2 - {\delta _2}}}}{\mathop{\rm cov}} erB\left( {x,{\rho ^{1/2 + {\delta _2}}}} \right)} {\int_{\Pi  \cap {B_{{}_{{\rho ^{1/2 - {\delta _2}}}}}}} {{{\left| {{e^{it{\rm H}}}f{\psi _2}\left( {\frac{t}{R}} \right)} \right|}^2}} }  \nonumber\\
  &\le \sum\limits_{B{}_{{\rho ^{1/2 - {\delta _2}}}}{\mathop{\rm cov}} erB\left( {x,{\rho ^{1/2 + {\delta _2}}}} \right)} {{R^{O\left( {{\delta _2}} \right)}}{{\left( {\frac{{{R^{1/2 - {\delta _2}}}}}{{{\rho ^{1/2 - {\delta _2}}}}}} \right)}^{ - 1}}\int_{\Pi  \cap {B_{{}_{{R^{1/2 - {\delta _2}}}}}}} {{{\left| {{e^{it{\rm H}}}f{\psi _2}\left( {\frac{t}{R}} \right)} \right|}^2}} }  \nonumber\\
  &\le {R^{O\left( {{\delta _2}} \right)}}{\left( {\frac{{{R^{1/2}}}}{{{\rho ^{1/2}}}}} \right)^{ - 1}}\int_{\Pi  \cap 2B} {{{\left| {{e^{it{\rm H}}}f{\psi _2}\left( {\frac{t}{R}} \right)} \right|}^2}}, \nonumber
\end{align}
here we used the fact that on $\Pi$ which passing through $B$, ${{e^{it{\rm H}}}f{\psi _2}\left( {\frac{t}{R}} \right)}$ is essentially supported in $\Pi \cap 2B$, see also \cite{G}. By \cite{G},  $\Pi \cap B \cap N_{\rho^{1/2+\delta_{2}}}(Z) \subset {N_{{\rho ^{ 1/2 + {\delta _2}}}}}\left( {\Pi  \cap Z} \right) \cap \Pi  \cap 2B,$ and ${N_{{\rho ^{ 1/2 + {\delta _2}}}}}\left( {\Pi  \cap Z} \right) \cap \Pi  \cap 2B$ can be covered by ${R^{O\left( {{\delta _2}} \right)}}$ balls $\Pi \cap B(x,\rho^{1/2+\delta_{2}}),x\in B \cap \Pi$, so we get the bound
\[\left\| {{e^{it{\rm H}}}f{\psi _2}\left( {\frac{t}{R}} \right)} \right\|_{{L^2}\left( {\Pi  \cap B \cap {N_{{\rho ^{1/2 + {\delta _2}}}}}\left( Z \right)} \right)}^2 \le {R^{O\left( {{\delta _2}} \right)}}{\left( {\frac{{{R^{1/2}}}}{{{\rho ^{1/2}}}}} \right)^{ - 1}}\left\| {{e^{it{\rm H}}}f{\psi _2}\left( {\frac{t}{R}} \right)} \right\|_{{L^2}\left( {\Pi  \cap 2B} \right)}^2,\]
(\ref{Eq72}) is obtained by integrating over all $\Pi $ paralleled to $\left( {{a_1}, {a_2},0} \right)$ and this completes the proof of Theorem \ref{Theorem7.1}.

In the proof of Theorem \ref{Theorem6.1}, we also used the following generalization of Lemma 3.4 in \cite{G}:
\begin{theorem}\label{Theorem7.2}
Suppose that $f$ is concentrated in a set of wave packets ${\rm T}$ and that for every $\left( {\theta ,\nu } \right) \in {\rm T},$ ${T_{\theta ,\nu }} \cap B\left( {z,r} \right) \ne \emptyset ,$ $z = \left( {{x_0},{t_0}} \right),{t_0} \le R,$ for some radius $r \sim {R^{1/2 + \delta_{2} }}.$ Then
\begin{equation}\label{Eq76}
\left\| {{e^{it{\rm H}}}f{\psi _2}\left( {\frac{t}{R}} \right)} \right\|_{{L^2}\left( {B\left( {z,10r} \right)} \right)}^2 \sim r\left\| f \right\|_{{L^2}}^2.
\end{equation}
\end{theorem}
\textbf{Proof:} Suppose $z = \left( {{x_0},{t_0}} \right),$ for each $t$ in the range ${t_0} - r \le t \le {t_0} + r,$ each $\left( {\theta ,\nu } \right) \in {\rm T},$ ${T_{\theta ,\nu }} \cap \left( {{\mathbb{R}^2} \times \left\{ t \right\}} \right) \subset B\left( {{x_0},5r} \right),$ therefore, (\ref{Eq76}) follows from the facts that
\begin{align}
 \left\| {{e^{it{\rm H}}}f{\psi _2}\left( {\frac{t}{R}} \right)} \right\|_{{L^2}\left( {B\left( {z,10r} \right)} \right)}^2 &\ge \left\| {{e^{it{\rm H}}}f{\psi _2}\left( {\frac{t}{R}} \right)} \right\|_{{L^2}\left( {B\left( {{x_0},5r} \right) \times \left[ {{t_0} - r,{t_0} + r} \right]} \right)}^2 \nonumber\\
 &= \int_{{t_0} - r}^{{t_0} + r} {\int_{B\left( {{x_0},5r} \right)} {{{\left| {{e^{it{\rm H}}}f} \right|}^2}dx} } {\left| {{\psi _2}\left( {\frac{t}{R}} \right)} \right|^2}dt \nonumber\\
  &= \int_{{t_0} - r}^{{t_0} + r} {\int_{{\mathbb{R}^2}} {{{\left| {{e^{it{\rm H}}}f} \right|}^2}dx} } {\left| {{\psi _2}\left( {\frac{t}{R}} \right)} \right|^2}dt = \left\| f \right\|_{{L^2}}^2\int_{{t_0} - r}^{{t_0} + r} {{{\left| {{\psi _2}\left( {\frac{t}{R}} \right)} \right|}^2}dt}  \nonumber\\
  &\ge \left( {r/2} \right)\left\| f \right\|_{{L^2}}^2, \nonumber
\end{align}
and
\begin{align}
 \left\| {{e^{it{\rm H}}}f{\psi _2}\left( {\frac{t}{R}} \right)} \right\|_{{L^2}\left( {B\left( {z,10r} \right)} \right)}^2 &\le \int_{{t_0} - 10r}^{{t_0} + 10r} {\int_{{\mathbb{R}^2}} {{{\left| {{e^{it{\rm H}}}f} \right|}^2}dx} } {\left| {{\psi _2}\left( {\frac{t}{R}} \right)} \right|^2}dt = \left\| f \right\|_{{L^2}}^2\int_{{t_0} - 10r}^{{t_0} + 10r} {{{\left| {{\psi _2}\left( {\frac{t}{R}} \right)} \right|}^2}dt}  \nonumber\\
  &\le \left\| f \right\|_{{L^2}}^2\int_{{t_0} - 10r}^{{t_0} + 10r} {1dt}  \le 20r\left\| f \right\|_{{L^2}}^2. \nonumber
\end{align}


\begin{thebibliography}{99}
\bibitem{C} Carleson L., Some analytic problems related to statistical mechanics. Euclidean harmonic analysis. Springer, Berlin, Heidelberg, 1980: 5-45.

\bibitem{CLV} Cho C. H., Lee S., Vargas A., Problems on pointwise convergence of solutions to the Schr\"{o}dinger equation. Journal of Fourier Analysis and Applications, 2012, 18(5): 972-994.

\bibitem{DK} Dahlberg B. E. J., Kenig C. E., A note on the almost everywhere behavior of solutions to the Schr\"{o}dinger equation. Harmonic analysis. Springer, Berlin, Heidelberg, 1982: 205-209.

\bibitem{DN} Ding Y., Niu Y., Weighted maximal estimates along curve associated with dispersive equations. Analysis and Applications, 2017, 15(02): 225-240.

\bibitem{DL} Du X., Li X., $ L^ p $-estimates of maximal function related to Schr\"{o}dinger Equation in $\mathbb {R}^ 2$. arXiv preprint arXiv:1508.05437, 2015.

\bibitem{DGL} Du X., Guth L., Li X., A sharp Schr\"{o}dinger maximal estimate in $\mathbb{R}^{2}$. Annals of Mathematics, 2017: 607-640.

\bibitem{G} Guth L., Restriction estimates using polynomial partitioning II, preprint, 2016. arXiv preprint arXiv:1603.04250.

\bibitem{SS} Sj\"{o}gren P., Sj\"{o}lin P., Convergence [properties for the time-dependent Schr\"{o}dinger equation. Annales Academire Scientiarurn Fennicre, Series A. I. Mathematica, 1987, 14: 13-25.

\bibitem{S} Sj\"{o}lin P., Regularity of solutions to the Schr\"{o}dinger equation. Duke Mathematical journal, 1987, 55(3): 699-715.

\bibitem{V} Vega L., Schr\"{o}dinger equations: pointwise convergence to the initial data. Proceedings of the American Mathematical Society, 1988, 102(4): 874-878.

\end{thebibliography}
\end{document}